\documentclass[11pt]{amsart}
\usepackage{graphicx}

\newtheorem{thm}{Theorem}[section]
\newtheorem{prop}[thm]{Proposition}
\newtheorem{lemma}[thm]{Lemma}
\newtheorem{cor}[thm]{Corollary}
\newtheorem{conj}[thm]{Conjecture}

\newtheorem*{appLemma}{Lemma}

\newcommand{\LL}[1]{\ensuremath{\mathcal{L}^{#1}}}
\newcommand{\Lt}[1]{\ensuremath{\mathcal{\widetilde{L}}^{#1}}}

\newcommand{\CC}[1]{\ensuremath{\mathbf{C}^{#1}}}
\newcommand{\CP}[1]{\ensuremath{\mathbf{CP}^{#1}}}
\newcommand{\PP}[1]{\ensuremath{\mathbf{P}^{#1}}}
\newcommand{\RP}[1]{\ensuremath{\mathbf{RP}^{#1}}}

\newcommand{\HH}{\ensuremath{\mathcal{H}}}
\newcommand{\Hh}{\ensuremath{\mathbf{H}}}

\newcommand{\Sym}[1]{\ensuremath{\mathcal{S}_{#1}}}
\newcommand{\G}[1]{\ensuremath{\mathcal{G}_{#1}}}
\newcommand{\Z}[1]{\ensuremath{\mathbf{Z}_{#1}}}

\newcommand{\arrGap}{5pt}

\newcommand{\basinWidth}{3.5in}

\newcommand{\dd}[1]{\frac{d}{d #1}}
\newcommand{\pdd}[2]{\frac{\partial{#1}}{\partial{#2}}}

\newcommand{\gT}{\ensuremath{\widetilde{g}}}

\allowdisplaybreaks

\begin{document}

\title{A family of critically finite maps with symmetry}
\author{Scott Crass}
\address{Mathematics Department\\
California State University, Long Beach\\
Long Beach, CA  90840-1001}
\date{\today}
\email{scrass@csulb.edu} 

\begin{abstract}

The symmetric group \Sym{n} acts as a reflection group on \CP{n-2}
(for $n\geq 3$) .  Associated with each of the $\binom{n}{2}$
transpositions in \Sym{n} is an involution on \CP{n-2} that
pointwise fixes a hyperplane---the mirrors of the action.  For
each such action, there is a unique \Sym{n}-symmetric holomorphic
map of degree $n+1$ whose critical set is precisely the collection
of hyperplanes. Since the map preserves each reflecting
hyperplane, the members of this family are critically-finite in a
very strong sense.  Considerations of symmetry and
critical-finiteness produce global dynamical results: each map's
Fatou set consists of a special finite set of superattracting
points whose basins are dense.

\end{abstract}
\maketitle

\section{Overview}

Complex dynamics in several dimensions has been the object of
considerable recent study.  Some specialized previous work in this
field treats a variety of maps that share a common property: they
respect the action of a finite group on a complex projective
space. (See \cite{sextic}, \cite{quintic}, \cite{octic}.)  The
nature of these investigations leads to a consideration of issues
pertaining to global dynamics. While the most significant
dynamical claims possess experimental support, they remain
theoretical conjectures. The current project stems from a desire
to find symmetrical maps with interesting geometry and tractable
dynamics.  Its first fruit is an infinite family of special maps
each of whose members respect the action of the symmetric group
\Sym{n}.  In fact, for each $n \geq 3$, there is a unique
holomorphic map $g$ on \CP{n-2} whose critical set consists of an
\Sym{n} orbit of $\binom{n}{2}$ hyperplanes that $g$ preserves.
This leads to a strong form of critical finiteness that yields
several global dynamical results of the type that eluded earlier
undertakings.

The treatment develops in three stages:

\begin{enumerate}

\item some background on special actions of \Sym{n} and their
associated symmetrical maps

\item proofs that the special family of critically-finite maps
with \Sym{n} symmetry exists and that each member is unique and
holomorphic

\item proofs of claims concerning the dynamics of the
maps (in the cases $n=3,4$)  Specifically, each member has a certain
attractor with dense basins.  When $n>4$, the claim concerning the
attractor is conjectured.

\end{enumerate}
Finally, some graphical results for low-dimensional cases appear.

\section{\Sym{n} acts on \CP{n-2}}

The permutation action of the symmetric group \Sym{n} on \CC{n}
preserves the hyperplane
$$\HH = \Biggl\{\sum_{k=1}^n x_k =0 \Biggr\} \simeq \CC{n-1}$$
and, thereby, restricts to a faithful $(n-1)$-dimensional
irreducible representation. This action on \CC{n-1} projects
one-to-one to a group \G{n} on $\Hh:=\PP{}\HH \simeq \CP{n-2}$.

\subsection{Special orbits and reflection hyperplanes}

The smallest \G{n} orbit consists of the $n$ points
$$[1-n,1,\dots,1],\dots,[1,\dots,1,1-n].$$
(Square brackets indicate points in projective space.)

Corresponding to the $\binom{n}{2}$ transpositions $(ij)$ in
\Sym{n} are $\binom{n}{2}$ involutions
$$ x_i \longleftrightarrow x_j$$
on \Hh\ that generate \G{n} as a \emph{complex reflection group}.
Each generating involution fixes the point
$$
[0,\dots,0, \overbrace{1}^i,0,\dots,0, \overbrace{-1}^j,0,\dots,0]
$$
and pointwise fixes the companion hyperplane $ \{x_i=x_j\}.$ This
point-hyperplane pair gives the only fixed points of the
involution. They form \G{n} orbits of size $\binom{n}{2}$. For
ease of reference, use the term ``$\binom{n}{2}$-hyperplane."

\subsection{Coordinates}

The transformation $A:\CC{n}\rightarrow \CC{n-1}$ given by
$$
u=Ax \qquad A= \begin{pmatrix}
1&0&\dots&0&-1\\
0&1&\dots&0&-1\\
\vdots&\vdots&\vdots&\vdots&\vdots\\
0&0&\dots&1&-1
\end{pmatrix}
=\begin{pmatrix}a_{ij}\end{pmatrix} \qquad a_{ij}=
\begin{cases}
1&i=j\\
-1&j=n\\
0&\text{otherwise}
\end{cases}
$$
gives a special system of $n-1$ coordinates on \Hh\ where the
$n$-point orbit is
$$[1,0,\dots,0],\dots,[0,\dots,0,1],[1,\dots,1].$$
Note that the null space of $A$ is the euclidean orthogonal
complement to \Hh.  This change of coordinates has an ``inverse"
$$
x=Bu \qquad B= \begin{pmatrix}
1-n&1&1&\dots&1&1\\
1&1-n&1&\dots&1&1\\
\vdots&\vdots&\vdots&\vdots&\vdots&\vdots\\
1&1&1&\dots&1&1-n\\
1&1&1&\dots&1&1
\end{pmatrix}
$$
which gives
$$AB=-n\,I_{n-1} \qquad BA=\mathbf{1}_n - n\,I_n$$
where $I_m$ is the $m\times m$ identity and $\mathbf{1}_n$ is the
$n\times n$ matrix each entry of which is $1$.  Accordingly, $A$
and $B$ induce isomorphisms between \Hh\ and \CP{n-2}.

In $u$-coordinates, the $\binom{n}{2}$-hyperplanes are the $n-1$
coordinate hyperplanes $\{u_k=0\}$ and the $\binom{n-1}{2}$ spaces
$\{u_k=u_\ell\}.$  The points determined by the intersections of
the $\binom{n}{2}$-hyperplanes play a central role in subsequent
developments. Their description is especially simple in $u$. (See
Table~\ref{tab:intpts2}.)  With one exception, each orbit consists
of points $p_k$ and $q_k$ with complementary coordinates.

\begin{table}[h]
$$
\begin{array}{c|c|c}

&\text{representative points}&\\
n&\text{on $n-2$ hyperplanes}&\text{orbit size}\\
\hline
&&\\
2\,m-1&p_k=[\overbrace{1,\dots,1}^{k},
\overbrace{0,\dots,0}^{n-k-1}]&\binom{n-1}{k}+\binom{n-1}{k-1}=\binom{n}{k}\\[\arrGap]
&q_k=[\overbrace{1,\dots,1}^{n-k},
\overbrace{0,\dots,0}^{k-1}]&k=1,\dots,m-1\\[10pt]
2\,m&p_k=[\overbrace{1,\dots,1}^{k},
\overbrace{0,\dots,0}^{n-k-1}]&\begin{cases}
\binom{n}{k}& k<m\\
\frac{1}{2}\binom{n}{k}=\binom{n-1}{k-1}& k=m
\end{cases}\\[\arrGap]
&q_k=[\overbrace{1,\dots,1}^{n-k},
\overbrace{0,\dots,0}^{k-1}]&k=1,\dots,m
\end{array}
$$

\caption{Points determined by intersections of
$\binom{n}{2}$-hyperplanes}

\label{tab:intpts2}

\end{table}

Relative to the $u$ space, \G{n} is generated over the permutation
action \G{n-1} of \Sym{n-1} on the $u_k$ by means of the
involution
$$T=\begin{pmatrix}
-1&0&0&\dots&0&0\\
-1&1&0&\dots&0&0\\
\vdots&\vdots&\vdots&\vdots&\vdots&\vdots\\
-1&0&0&\dots&1&0\\
-1&0&0&\dots&0&1
\end{pmatrix}
$$
that transposes the pair $\{p_1,q_1\}$ and fixes the remaining
members of the $n$-point orbit.  Note that $T$ is the $u$ version
of the transformation
$$x_1 \longleftrightarrow x_n.$$

\section{\G{n} equivariants}

Consider a map
$$f=[f_1,\dots,f_{n-1}]$$
from \Hh\ to itself given by homogeneous polynomials in
$$u=(u_1,\dots,u_{n-1})$$
of degree $r$.  In general, $f$ can be meromorphic; that is, for
some $p \in \CC{n-1}$, $f(p)=0$ for every lift of $f$ to \CC{n-1}.
We say that $f$ is \G{n}-\emph{equivariant} when it sends a group
orbit to a group orbit.  Algebraically, this means that $f$
commutes with every element of \G{n}. Obviously, $f$ is
\G{n-1}-equivariant as well. It readily follows that each
component $f_k$ is invariant under the stabilizer $Z_k$ of $u_k$.
Thus, we can express a component by
$$f_k=\sum_{\ell=0}^r u_k^{r-\ell} A_{k,\ell}$$
where $A_{k,\ell}$ is a degree-$\ell$ $Z_k$ invariant.
Accordingly, each $A_{k,\ell}$ is taken to be a polynomial in the
elementary symmetric functions in the complementary variables
$$\widehat{u}_k = (u_1, \dots,u_{k-1},u_{k+1},\dots u_{n-1}).$$
Alternatively, we can employ the elementary symmetric functions in
$u$ when expressing $A_{k,\ell}$. This is a matter of expressing a
polynomial in $\widehat{u}_k$ in terms of a polynomial in $u$ and
a polynomial in $\widehat{u}_k$ with lower degree. Specifically,
let $\widehat{S}_m$ and $S_m$ be the degree-$m$ elementary
symmetric functions in $\widehat{u}_k$ and $u$ respectively.
Taking $S_0=1$, the relations
$$\widehat{S}_m = S_m - u_k\,\widehat{S}_{m-1}$$
give a reductive scheme for the replacement process.

An immediate consequence of \G{n-1} equivariance is that
$$A_{j,\ell}=A_{k,\ell}:=A_\ell \qquad \text{for all}\ j,k, \ell.$$
We can say a bit more concerning the form that \G{n} equivariants
take.

First, consider a point $a$ that some element $M \in \G{n}$ fixes
and observe that
$$Mf(a)=f(Ma)=f(a).$$
Thus, $f$ either sends $a$ to another fixed point of $M$ or blows
up at $a$---that is, for any lift $\widetilde{f}$ and
$\widetilde{a}$ of $f$ and $a$ to \CC{n-1},
$\widetilde{f}(\widetilde{a})=0$. Applying this condition to the
$\binom{n}{2}$-hyperplanes, provided that $a$ is not a point of
indeterminacy, each point on such a hyperplane must map to a point
that is fixed by the involution that fixes the hyperplane
pointwise. The only place for the image of such a point is on the
hyperplane itself or its companion point. Under a holomorphic map,
the image cannot be the companion point---this would force the
entire hyperplane to collapse to the point.  So, a holomorphic
\G{n} equivariant $f$ sends an $\binom{n}{2}$-hyperplane to
itself. This circumstance forces $f_k$ to be divisible by $u_k$
and, thereby, requires the terms $A_r$ to be a power of $S_{n-1}$
or to vanish. In particular, when $r\leq n-1$, $A_r=0$ so that
$$f_k=u_k \sum_{\ell=0}^{r-1} u_k^{r-\ell-1} A_\ell.$$

By design, the map $f$ has \G{n-1} symmetry.  To be fully
 \G{n}-equivariant, the map must commute with $T$ as well.  This
condition places strong restrictions on the $A_\ell.$  The general
form they take might be an interesting result, but not one taken
up by the current investigation.  Here, the quest is for a family
of \G{n} equivariants with very special properties.

\section{Reflection hyperplanes as critical sets: existence,
uniqueness, and holomorphy}

Explicit computation in low-degree cases reveals the existence of
a unique holomorphic \G{n} equivariant whose critical set is
precisely the $\binom{n}{2}$-hyperplanes counted with multiplicity
two. These maps conform to a general formula.  Let
$$g=[g_1, \dots, g_{n-1}]$$
where
$$
g_\ell = u_\ell^3\,G_\ell,\qquad G_\ell = \sum_{k=0}^{n-2}
(-1)^k\,\frac{k+1}{k+3}\,u_\ell^k\,S_{n,n-2-k},
$$
and $S_{n,\ell}$ is the degree-$\ell$ elementary symmetric
function in $u_1,\dots,u_{n-1}$. In the degree-$0$ case, take
$S_{n,0}=1.$ By construction, each $g$ is equivariant under the
group \G{n-1} that permutes the $u_\ell$. In addition, the
$u_\ell^3$ factor in each coordinate implies that the maps are
doubly critical on $n-1$ of the
$\binom{n}{2}$-hyperplanes---namely, where $u_\ell=0$. Were $g$ to
commute with the transformation $T$ that generates \G{n} over
\G{n-1}, symmetry would provide for double criticality on the
remaining $\binom{n-1}{2}$ of the
$\binom{n}{2}$-hyperplanes---where $u_j=u_k.$ Moreover, since a
degree-(n+1) map in $n-1$ variables has a critical set whose
degree is
$$(n-1)\,n=2\,\binom{n}{2},$$
$g$'s critical set would consist exclusively of the
$\binom{n}{2}$-hyperplanes.

This section develops rather technical arguments for three main
results. According to Theorem~\ref{thm:gEqv}, the
$\binom{n}{2}$-hyperplanes form $g$'s critical set with
multiplicity two.  Moreover, Theorem~\ref{thm:gUnique} informs us
that there is only one such map for each \G{n} action.
Theorem~\ref{thm:gHlm} states that each $g$ is holomorphic on \Hh\
which implies that $g$ preserves each $\binom{n}{2}$-hyperplane
\LL{} rather than collapse \LL{} to a lower-dimensional variety; a
contraction would force the map to blow up.

Thus, $g$ is a family of maps each member of which is holomorphic,
doubly-critical on the $\binom{n}{2}$-hyperplanes, and
critically-finite.  As a standing assumption, let $n \geq 3$.

\begin{thm} \label{thm:gEqv}

The respective $g$ is $T$-equivariant, hence, \G{n}-equivariant.

\end{thm}

\begin{thm} \label{thm:gUnique}

Under the action of \G{n}, $g$ is the unique rational map of
degree $n+1$ for which each $\binom{n}{2}$-hyperplane is doubly
critical.

\end{thm}

\begin{thm} \label{thm:gHlm}

Each member of the family $g$ is holomorphic on \Hh.

\end{thm}

\begin{proof}[Proof of Theorem~\ref{thm:gEqv}]

Propositions~\ref{prop:g1oT} and \ref{prop:g2oT} below establish
that $g$ is symmetric under $T$ as well as under \G{n-1}. Since
$T$ generates \G{n} over \G{n-1}, $g$ is \G{n}-equivariant.

\end{proof}

The proofs of the propositions rely on a formula that describes
how the elementary symmetric functions transform under $T$.  This
result was found by pattern detection in low-degree cases.  For
simplicity of appearance, express the functions $S_{n,k}(u)$ in
the suppressed form $S_{n,k}$.

\begin{lemma} \label{lm:Snk}

For $k\leq n$, the \G{n-1} invariants $S_{n,k}$ transform under
$T$ according to
$$
S_{n,k}(T u) = \sum_{\ell=0}^k (-1)^\ell \binom{n-k+\ell}{n-k}
u_1^\ell\,S_{n,k-\ell}.
$$

\end{lemma}

\begin{proof}

Proofs of several technical lemmas appear in the appendix.

\end{proof}

The argument for the $T$-equivariance of $g$ examines the
coordinates individually.

\begin{prop} \label{prop:g1oT}

The factor $G_1$ of $g_1$ is $T$-invariant (in a linear as well as
projective sense).

\end{prop}

\begin{proof}

The proof amounts to manipulation of sums.  Since $n$ is fixed
here, let $S_{k}=S_{n,k}$.  Consider

\begin{equation*}
G_1(T u) = \sum_{k=0}^{n-2} (-1)^k\,\frac{k+1}{k+3}\,
  (-u_1)^k\,S_{n-2-k}(T u) =
  \sum_{k=0}^{n-2} \frac{k+1}{k+3}\,
  u_1^k\,S_{n-2-k}(T u).
\end{equation*}
By Lemma~\ref{lm:Snk},
\begin{align*}
G_1(T u) =&\ \sum_{k=0}^{n-2} \frac{k+1}{k+3}\,
  u_1^k\Biggl(
  \sum_{\ell=0}^{n-2-k} (-1)^\ell
  \binom{n-(n-2-k)+\ell}{n-(n-2-k)} u_1^\ell\,S_{n-2-k-\ell}
  \Biggr)\\
=&\ \sum_{k=0}^{n-2} \frac{k+1}{k+3}
  \Biggl(
  \sum_{\ell=0}^{n-2-k} (-1)^\ell
  \binom{k+\ell+2}{k+2} u_1^{k+\ell}\,S_{n-2-(k+\ell)}
  \Biggr).
\end{align*}
Setting $m=k+\ell$,
\begin{align*}
G_1(T u)
=&\ \sum_{k=0}^{n-2} \frac{k+1}{k+3}
  \Biggl(
  \sum_{m=k}^{n-2} (-1)^{m-k}\,
  \binom{m+2}{k+2}\,u_1^m\,S_{n-2-m}\Biggr).
\end{align*}
Reversing the order of summation,
\begin{align*}
G_1(T u)
=&\ \sum_{m=0}^{n-2}
  \Biggl(
  \sum_{k=0}^m (-1)^{m-k}\,\frac{k+1}{k+3}
  \binom{m+2}{k+2}\Biggr)\,u_1^m\,S_{n-2-m}\\
=&\ \sum_{m=0}^{n-2} (-1)^m\,(m+2)!
  \Biggl(
  \sum_{k=0}^m (-1)^k
  \frac{k+1}{(k+3)!} \frac{1}{(m-k)!}\Biggr) \,u_1^m\,S_{n-2-m}.
\end{align*}
Lemma~\ref{lm:specSum} below gives the sum over $k$:
\begin{align*}
G_1(T u)
=&\ \sum_{m=0}^{n-2} (-1)^m\,(m+2)!\,
  \frac{m+1}{(m+3)!}\,u_1^m\,S_{n-2-m}\\
=&\ \sum_{m=0}^{n-2} (-1)^m\,
  \frac{m+1}{m+3}\,u_1^m\,S_{n-2-m}\\
=&\ G_1(u).
\end{align*}

\end{proof}

\begin{lemma} \label{lm:specSum}

$$
\sum_{k=0}^m (-1)^k\,\frac{k+1}{(k+3)!\,(m-k)!} =
\frac{m+1}{(m+3)!}.
$$

\end{lemma}

\begin{proof}

See the appendix.

\end{proof}

\begin{cor}

Each $g$ is $T$-equivariant in the first coordinate.

\end{cor}

\begin{proof}

Let $[\cdot]_1$ specify a map's first coordinate.  Then
$$g_1\circ T = -u_1^3\,G_1\circ T = -u_1^3\,G_1 =
[T\circ g]_1.$$

\end{proof}

To establish overall $T$-equivariance, it suffices to consider the
behavior of $g$ under $T$ in just the second coordinate.  This
follows directly from the commutativity of $T$ and the members
$\tau_{2,m} \in \G{n}$ that simply transpose the second and $m$th
basis elements:
$$
[0,1,0,\dots,0]\stackrel{\tau_{2,m}}{\longleftrightarrow}
[0,0,\dots,0,\underbrace{1}_m,0,\dots,0]
$$
provided that $m\neq 1, n$.  Expressed in terms of \Sym{n}, this
amounts to the commutativity of the disjoint transpositions $(1n)$
and $(2m)$.  So, noting that $g$ is \G{n-1}-equivariant, hence,
$\tau_{2,m}$-equivariant, and given that $g$ is $T$-equivariant in
its second coordinate,
\begin{align*}
g_m \circ T =&\ [g\circ T]_m \\
=&\ [(\tau_{2,m}\circ g \circ \tau_{2,m})\circ T]_m =\
[\tau_{2,m}\circ (g \circ T \circ \tau_{2,m})]_m =\ [g \circ T
\circ \tau_{2,m}]_2\\
=&\ [g\circ T]_2 \circ \tau_{2,m} =\ [T\circ g]_2 \circ \tau_{2,m}
=\ [T\circ g\circ \tau_{2,m}]_2\\
=&\ [\tau_{2,m}\circ T\circ g]_2 =\ [T\circ g]_m.
\end{align*}

\begin{prop} \label{prop:g2oT}

The second coordinate of $g$ satisfies the equivariance condition
$$g_2 \circ T = [T \circ g]_2.$$

\end{prop}

\begin{proof}

First, express $g_2 \circ T$ in a way that's useful for comparison
to $[T \circ g]_2$.  Again, set $S_k=S_{n,k}$. Applying
Lemma~\ref{lm:Snk},
\begin{align*}
g_2(T u) =&\ (u_2-u_1)^3\sum_{k=0}^{n-2} (-1)^k\,\frac{k+1}{k+3}\,
  (u_2-u_1)^k\,S_{n-2-k}(T u) \\
=&\ \sum_{k=0}^{n-2} \frac{k+1}{k+3}\,
  \Biggl( \sum_{\ell=0}^{n-2-k} (-1)^{k+\ell}\,\binom{k+2+\ell}{k+2}\,
     u_1^\ell\,S_{n-2-k-\ell}\Biggr)\,(u_2-u_1)^{k+3}.
\end{align*}
Setting $m=k+\ell$,
\begin{align*}
g_2(T u) =&\ \sum_{k=0}^{n-2} \frac{k+1}{k+3}\,
  \Biggl( \sum_{m=k}^{n-2} (-1)^m\,\binom{m+2}{k+2}\,
     u_1^{m-k}\,S_{n-2-m}\Biggr)\,(u_2-u_1)^{k+3}.
\end{align*}
Reversing the order of summation,
\begin{equation*}
g_2(T u) = \sum_{m=0}^{n-2} (-1)^m
  \Biggl( \sum_{k=0}^m \frac{k+1}{k+3}\,
    \binom{m+2}{k+2}\,
    u_1^{m-k}\,(u_2-u_1)^{k+3} \Biggr)\,\,S_{n-2-m}.
\end{equation*}
Lemma~\ref{lm:g2specSum} below establishes a useful identity for
the sum over $k$ so that
\begin{align*}
g_2(T u) =&\ u_2^3\sum_{m=0}^{n-2} (-1)^m\,
  \frac{m+1}{m+3}\,u_2^m\,S_{n-2-m} -
  u_1^3\sum_{m=0}^{n-2} (-1)^m\,
 \frac{m+1}{m+3}\,u_1^m\,S_{n-2-m}\\
& -  u_1\,u_2\,\sum_{m=0}^{n-2} (-1)^m\,( u_2^{m+1}-u_1^{m+1})\,
S_{n-2-m}.
\end{align*}
The first two terms are $g_2(u)$ and $g_1(u)$ respectively. Since
their difference amounts to $[Tg(u)]_2$,
\begin{align*}
\bigl[T g(u)\bigr]_2 - g_2(T u) =&\ u_1\,u_2\,
  \sum_{m=0}^{n-2} (-1)^m
  (u_2^{m+1}-u_1^{m+1})\,S_{n-2-m}.
\end{align*}
Adding and subtracting $-u_1\,u_2\,S_{n-1}$ on the right,
\begin{align*}
\bigl[T g(u)\bigr]_2 - g_2(T u) =&\ u_1\,u_2\,\Bigl(
  \bigl(-S_{n-1}
    + \sum_{m=0}^{n-2}(-1)^m\,u_2^{m+1}\,S_{n-2-m}\bigr) \\
&\ - \bigl(-S_{n-1}
    + \sum_{m=0}^{n-2}(-1)^m\,u_1^{m+1}\,S_{n-2-m}\bigr)
\Bigr).
\end{align*}
Let $m=p-1$, while, for the apparent variables $u_1$ and $u_2$,
set $x=u_2$ and $y=u_1$.  The result is
\begin{align*}
\bigl[T g(u)\bigr]_2 - g_2(T u) =&\ x\,y\,\Bigl(
  -\sum_{p=0}^{n-1}(-1)^p\,x^p\,S_{n-1-p}
   + \sum_{p=0}^{n-1}(-1)^p\,y^p\,S_{n-1-p}
\Bigr) \\
=&\ x\,y\,\Bigl(
 -\prod_{k=1}^{n-1} (u_k-x) + \prod_{k=1}^{n-1} (u_k-y)
\Bigr).
\end{align*}
Thus, when $x=u_2$ and $y=u_1$,
$$\bigl[T g(u)\bigr]_2 - g_2(T u) = 0.$$

\end{proof}

\begin{lemma} \label{lm:g2specSum}

\begin{equation*}
\sum_{k=0}^m \frac{k+1}{k+3}\,\binom{m+2}{k+2}\,
  u_1^{m-k}\,(u_2-u_1)^{k+3}
= \frac{m+1}{m+3}\,\bigl(u_2^{m+3}-u_1^{m+3}\bigr)-
  u_1\,u_2\,\bigl(u_2^{m+1}-u_1^{m+1}\bigr).
\end{equation*}

\end{lemma}

\begin{proof}

See the appendix.

\end{proof}

Now we turn to the matter of uniqueness.

\begin{proof}[Proof of Theorem~\ref{thm:gUnique}]

Suppose that
$$h=[h_1,\dots,h_{n-1}]$$
is a map of this type. The strategy is to compare $g$ to $h$ in
terms of $u$ coordinates.  Since $h$ is \G{n-1}-equivariant and
doubly critical on each $\{u_k=0\}$, the components of $h$ have
the form
$$h_k=u_k^3\,H_k.$$
Furthermore, each $H_k$ is a degree-$(n-2)$ invariant under an
\Sym{n-2}-isomorphic subgroup of \G{n-1}, namely, the stabilizer
of $u_k$.  It follows that we can express these polynomials by
$$H_k=\sum_{\ell=0}^{n-2} u_k^{n-2-\ell}\,V_\ell$$
where $V_\ell$ is a \G{n-1} invariant of degree $\ell$.

By \G{n-1} symmetry, we can examine a single component: $h_1$,
say.  Now, consider $V_{n-2}$.  In the event that $u_1$ divides
$V_{n-2}$, the associated component takes the form
$$h_1=u_1^4\,\widehat{H_1}.$$
  But this implies that $\{u_1=0\}$ is triply critical which
is at odds with the assumption that $h$ is doubly critical on the
$\binom{n}{2}$-hyperplanes.  By degree counting, the latter state
of affairs completely accounts for the critical set.

Accordingly, assume that $V_{n-2} \not\equiv 0$ when $u_1=0$. We
can now say that
$$V_{n-2}=u_1\,X + Y$$
where no monomial in $Y$ contains $u_1$. Hence, $Y$ is invariant
under the stabilizer in \G{n-1} of $u_1$. Lemma~\ref{lm:Y_1} below
reveals that $Y$ is divisible by each $u_k$ except $u_1$, of
course. Since the degree of $Y$ is $n-2$, this result implies that
$$Y = \alpha\prod_{k=2}^{n-1} u_k$$
where $\alpha \in \CC{}-\{0\}$. The \G{n-1} invariance of
$V_{n-2}$ requires that every element in the \G{n-1} orbit of $Y$
appears in $V_{n-2}$ and only these terms appear. Thus,
$$V_{n-2} = \alpha\,S_{n-2}.$$

Recalling the form of $g$, lift $g$ and $h$ to maps $\gT$ and
$\widetilde{h}$ on \CC{n-1} so that
$$G_1|_{u_1=0}=H_1|_{u_1=0}.$$
Also, we can lift \G{n-1} trivially to a linear group
$\widetilde{\G{}}$. Consequently, the $\widetilde{\G{}}$
equivariant $\gT-\widetilde{h}$ is either the zero map or is both
doubly critical along the $\binom{n}{2}$-hyperplanes and, as in
the case considered above, has the contrary property that its
first component is divisible by $u_1^4$.  Hence, the former case
is the only possibility so that $\widetilde{h}=\gT$.

\end{proof}

Evidently, $g$'s uniqueness is due to its \emph{full} \G{n}
symmetry---that is, to its $T$-equivariance in addition to its
symmetry under \G{n-1}. The proof of the following lemma makes
this explicit.

\begin{lemma} \label{lm:Y_1}

Define $Y$ as above.  For $k\neq 1$, $Y|_{u_k=0}=0$.

\end{lemma}

\begin{proof}

Let $k \neq 1$.   Equivariance under $T$ requires the components
of $h$ to satisfy the following identities:
\begin{align}
H_1 \circ T =&\ H_1 \label{eq:h1Cond} \\
(u_k-u_1)^3\,H_k \circ T =&\ u_k^3\,H_k - u_1^3\,H_1.
\label{eq:hkCond}
\end{align}
(To lessen clutter, suppress explicit mention of the variable $u$,
where possible.)  By (\ref{eq:h1Cond}),
$$
\sum_{\ell=0}^{n-2}(-u_1)^{n-2-\ell}\,V_\ell \circ T = H_1 \circ T
=H_1= \sum_{\ell=0}^{n-2} u_1^{n-2-\ell}\,V_\ell.
$$
From this we obtain
$$
V_{n-2} \circ T=V_{n-2} + \sum_{\ell=0}^{n-3}
\bigl(V_\ell-(-1)^\ell\,V_\ell \circ T \bigr)\,u_1^{n-2-\ell}
$$
which we can abbreviate to
\begin{equation} \label{eq:Vn-2}
V_{n-2} \circ T = V_{n-2} + u_1\,W_{n-3}.
\end{equation}

Turning to (\ref{eq:hkCond}),
\begin{align*}
(u_k-u_1)^3 \sum_{\ell=0}^{n-2}(u_k-u_1)^{n-2-\ell}\,V_\ell \circ
T =&\
\sum_{\ell=0}^{n-2}
\bigl(u_k^{n+1-\ell}-u_1^{n+1-\ell}\bigr)\,V_\ell \\
(u_k-u_1)^3\,V_{n-2}\circ T - (u_k^3-u_1^3)\,V_{n-2} =&\
\sum_{\ell=0}^{n-3} \Bigl(
\bigl(u_k^{n+1-\ell}-u_1^{n+1-\ell})\,V_\ell -
(u_k-u_1)^{n+1-\ell}\,V_\ell \circ T \Bigr). \\
\end{align*}
Expanding the first binomial on the left, using (\ref{eq:Vn-2}),
and rearranging gives
\begin{align*}
3\,u_1\,u_k\,(u_k-u_1)\,V_{n-2} =&\ (u_k-u_1)^3\,u_1\,W_{n-3} \\
& - \sum_{\ell=0}^{n-3} \Bigl(
\bigl(u_k^{n+1-\ell}-u_1^{n+1-\ell}\bigr)\,V_\ell -
(u_k-u_1)^{n+1-\ell}\,V_\ell \circ T \Bigr).
\end{align*}
Dividing through by the common factor $u_k-u_1$,
\begin{align*}
3\,u_1\,u_k\,V_{n-2} =&\ (u_k-u_1)^2\,u_1\,W_{n-3} \\
& - \sum_{\ell=0}^{n-3} \Biggl( \Biggl( \sum_{m=0}^{n-\ell}
u_1^m\,u_k^{n-\ell-m} \Biggr)\,V_\ell - (u_k-u_1)^{n-\ell}\,V_\ell
\circ T \Biggr).
\end{align*}
Restricting to $\{u_1=u_k\}$,
\begin{align*}
3\,u_1^2\,(V_{n-2}|_{u_1=u_k}) =&\ \sum_{\ell=0}^{n-3}
(n-\ell+1)\,u_1^{n-\ell}\,(V_\ell|_{u_1=u_k}) \\
=&\ u_1^3\,\sum_{\ell=0}^{n-3}
(n-\ell+1)\,u_1^{n-\ell-3}\,(V_\ell|_{u_1=u_k}).
\end{align*}
Note that this expression makes sense since $n\geq 3$.  Thus,
\begin{align*}
3\,(V_{n-2}|_{u_1=u_k}) =&\ u_1\,\sum_{\ell=0}^{n-3}
(n-\ell+1)\,u_1^{n-3-\ell}\,(V_\ell|_{u_1=u_k}).
\end{align*}
Finally, since $Y=V_{n-2}|_{u_1=0}$,
$$
Y|_{u_k=0}= \bigl(V_{n-2}|_{u_1=0}\bigr)\bigr|_{u_1=u_k} =
\bigl(V_{n-2}|_{u_1=u_k}\bigr)\bigr|_{u_1=0} = 0.
$$

\end{proof}

The upcoming proof of Theorem~\ref{thm:gHlm} exploits a
dimension-reducing process of restricting $g$ to intersections of
$\binom{n}{2}$-hyperplanes.  This cascade of intersections leads
to the special point-orbits determined by the hyperplanes.  At
these points, the map's behavior is explicitly computable.

\begin{proof}[Proof of Theorem~\ref{thm:gHlm}]

When $n=3$, $g$ is one-dimensional and hence, holomorphic.  As for
the non-trivial cases $n>3$, choose the ``literal" lift of $g$ to
\CC{n-1}:
$$\gT=(g_1, \dots, g_{n-1})$$
where
$$
g_\ell = u_\ell^3\,G_\ell \qquad \text{and} \qquad G_\ell =
\sum_{k=0}^{n-2} (-1)^k\,\frac{k+1}{k+3}\,
  u_\ell^k\,S_{n,n-2-k}.
$$

Let $X$ denote the union of the $\binom{n}{2}$-hyperplanes lifted
to hyperspaces through $0$ in \CC{n-1}. Suppose there is a point
$a \in \CC{n-1}$ where $\gT(a)=0$. By homogeneity,
$$(n+1)\,D\gT(u) = D\gT(u)\,u$$
where
$$
D\gT(u)=
\begin{pmatrix}
\frac{\partial{g_i}(u)}{\partial{g_j}(u)}
\end{pmatrix}
$$
is the Jacobian matrix of $\gT$. Thus, $\gT$ is critical at $a$.
 That is, $a$ is a zero eigenvector for $D\gT(a)$.  In this case,
the map collapses in the ``radial" direction defined by $a$. Since
$\gT$ is critical only on $X$, $a$ lies on one of $X$'s
constituent hyperspaces; call this hyperspace \LL{n-2}{} ($\simeq
\CC{n-2}$) and consider the restriction $\gT_{n-2}$ of $\gT$ to
\LL{n-2}{}. (Note that the action of \G{n} restricted to \LL{n-2}
is isomorphic to \Sym{n-2} so that $\gT_{n-2}$ is \emph{not} the
member of the family $g$ for dimension $n-2$ where the action is
that of \Sym{n-1}.)

Since
$$\gT_{n-2}(a)=0,$$
$a$ is a zero eigenvector for $D\gT_{n-2}(a)$; the critical set of
$\gT_{n-2}$ contains $a$. But, a zero eigenvector $v$ for
$D\gT_{n-2}(a)$ corresponds to a radial collapse in the $v$
direction so that $v$ is also a zero eigenvector for $D\gT(a)$.
But, as Lemma~\ref{lem:detVan} below describes,
$\det{D\gT_{n-2}(a)}$ does not vanish identically on \LL{n-2} so
that the critical set of $\gT_{n-2}$ is a proper algebraic subset
of $X$ and \LL{n-2}. Hence, the only possible location for $a$ is
where some hyperspace in $X$ different from \LL{n-2} intersects
\LL{n-2}. Denote this intersection by \LL{n-3}.

Further reducing the dimension, let
$$\gT_{n-3}=\gT|_{\LL{n-3}{}}$$
so that $\gT_{n-3}(a)=0$ and $a$ is critical for $\gT_{n-3}$.  As
above, $a$ belongs to the intersection of \LL{n-3} with a
hyperspace in $X$ that does not contain \LL{n-3}.

This reduction continues with the outcome at each stage that $a$
belongs to the intersection of $\binom{n}{2}$-hyperplanes.  When
the procedure arrives at dimension three, $a$ lies on two planes
through $0$ in \CC{n-1}---that is, a point in \CP{n-2}---that are
intersections of $\binom{n}{2}$-hyperspaces.  But,
Lemma~\ref{lem:g_pm} below implies that $\gT \neq 0$ at these
points.

\end{proof}

\begin{lemma} \label{lem:detVan}

For the restriction $\widehat{g}$ of $\gT$ to any space \LL{m} of
dimension $m \neq 0$ determined by the intersection of hyperspaces
in $X$, $\det{D\widehat{g}} \not\equiv 0$.

\end{lemma}

\begin{proof}

By the permutation action of \G{n-1} on the $u_k$, we can take
$$
\LL{m}=\Biggl(\bigcap_{k=1}^p \{u_k =0\}\Biggr) \bigcap
\Biggl(\bigcap_{\substack{i,j=1,\dots,n-m-p-2\\ p< \ell_i <
\ell_j}} \{u_{\ell_i} = u_{\ell_j}\}\Biggr).
$$
Any \LL{m} space that is partially determined by the intersection of
$p$ sets of the form $\{u_k=0\}$ belongs to the \G{n-1}-orbit of the
set specified above. Relabel the coordinates on \LL{m} so that the
restriction is expressed
$$
\widehat{g}(\widehat{u})=\gT_{m}|_{\LL{m}}(\widehat{u})=
\begin{pmatrix}
\widehat{g}_1(\widehat{u})\\
\vdots\\
\widehat{g}_m(\widehat{u})
\end{pmatrix}
=
\begin{pmatrix}
u_1^3\,\widehat{G}_1(\widehat{u})\\
\vdots\\
u_m^3\,\widehat{G}_m(\widehat{u})
\end{pmatrix}
\quad \text{where}\ \widehat{u}=\begin{pmatrix}
u_1\\
\vdots\\
u_{m}
\end{pmatrix}.
$$
Let $D\widehat{g}_i$ be the Jacobian of $\widehat{g}_i$ so that
$$
D\widehat{g}(\widehat{u}) =
\begin{pmatrix} D\widehat{g}_1(\widehat{u})\\
\vdots\\
D\widehat{g}_m(\widehat{u})
\end{pmatrix}.
$$
In order for $\det{D\widehat{g}}\equiv 0$, the set
$\{D\widehat{g}_i,\ i=1,\dots,m\}$ must be linearly dependent in
functional terms. To establish linear independence, consider the
relation
$$\sum_j^m a_j\,D\widehat{g}_j =0.$$

By homogeneity,
$$
(n+1)\,\widehat{g}_j(\widehat{u}) =
D\widehat{g}_j(\widehat{u})\,\widehat{u}
$$
and
$$
\lambda(\widehat{u}):=\sum_j^m a_j\,u_j^3\,\widehat{G}_j =
\sum_j^m a_j\,\widehat{g}_j = 0.
$$

But, on \LL{m} there are $m$ members of the \G{n-1} orbit of
$$p_1=(1,0,\dots,0),$$
namely,
$$\widehat{p}_k=(0,\dots,0,\overbrace{1}^k,0,\dots,0).$$
Since
$$\lambda(\widehat{p}_k)=a_k\,\widehat{G}_k(\widehat{p}_k),$$
the proof of Lemma~\ref{lem:g_pm} yields $a_k=0.$

\end{proof}

\begin{lemma} \label{lem:g_pm}

For the points $p_m$ that represent the orbits determined by the
intersections of $\binom{n}{2}$-hyperplanes,
$$\gT(p_m)\neq 0.$$

\end{lemma}

\begin{proof}

Recall that
$$
p_m=(\underbrace{1,\dots,1}_m,0,\dots,0) \qquad m=1,
\dots,\biggl\lceil \frac{n-1}{2}\biggr\rceil.
$$
It suffices to compute $G_1(p_m)$.

A straightforward calculation gives
$$
S_{n,k}(p_m) = \begin{cases} 0&k>m\\\binom{m}{k}&k\leq m
\end{cases}.
$$
With this,
\begin{align*}
G_1(p_m)=&\ \sum_{k=0}^{n-2} (-1)^k\,\frac{k+1}{k+3}\,
  S_{n,n-2-k}(p_m) \\
=&\ \sum_{k=n-2-m}^{n-2} (-1)^k\,\frac{k+1}{k+3}\,
\binom{m}{n-2-k}.
\end{align*}
Setting $p=n-2-k$,
\begin{align*}
G_1(p_m)=&\ \sum_{p=0}^{m} (-1)^{n-2-p}\,\frac{n-p-1}{n-p+1}\,
\binom{m}{p} \\
&=\ (-1)^n \sum_{p=0}^{m} (-1)^p\,\frac{n-p-1}{n-p+1}\,
\binom{m}{p}.
\end{align*}
From Lemma~\ref{lm:G1pm} below,
$$G_1(p_m)=(-1)^n\,\frac{2\,(-1)^{m-1}}{(n+1)\binom{n}{m}}\neq 0.$$

\end{proof}

\begin{lemma} \label{lm:G1pm}

$$
\sum_{p=0}^{m} (-1)^p\,\frac{n-p-1}{n-p+1}\, \binom{m}{p}
=\frac{2\,(-1)^{m-1}}{(n+1)\binom{n}{m}}.
$$

\end{lemma}

\begin{proof}

See the appendix.

\end{proof}

\section{Reflection hyperplanes as critical sets: Global dynamics}

Let \LL{n-3} generically denote an $\binom{n}{2}$-hyperplane and
let $X$ refer to the union of the \LL{n-3}. Where $m$ of the
\LL{n-3} intersect to form a \CP{n-2-m}, call the resulting space
\LL{n-2-m}. (Note that more than $m$ of the \LL{n-3} can pass
through an \LL{n-2-m}.)

Not only is $g$ critically-finite on $\Hh \simeq \CP{n-2}$ with
critical set consisting of the \LL{n-3} hyperplanes, the
restriction $g|_{\LL{n-2-m}}$ is also critically-finite, having a
collection of the \LL{n-3-m} for its critical set. In \cite{fs1},
such behavior is called \emph{strict} critical finiteness (Section
7).  In fact, \emph{all} of the \LL{n-3-m} on an \LL{n-2-m} are
critical for $g|_{\LL{n-2-m}}$ though not with the same
multiplicity.

\subsection{The Fatou set of $g$}

Following standard practice, the \emph{Fatou set} $F_g$ is where
the family of iterates $\{g^k\}$ is normal and the \emph{Julia
set} $J_g$ is the complement of $F_g$.

The behavior of $g$ on an \LL{n-3} plays a central dynamical role.
Again, lift $g$ to \CC{n-1}:
$$\widetilde{g}=(g_1, \dots, g_{n-1})$$
with
$$
g_\ell = u_\ell^3\,G_\ell \qquad \text{and} \qquad G_\ell =
\sum_{k=0}^{n-2} (-1)^k\,\frac{k+1}{k+3}\,
  u_\ell^k\,S_{n,n-2-k}.
$$
For a space $\LL{m}\subset \CP{k}$ lifted to \CC{k+1}, call the
lifted space \Lt{m+1}.

\begin{prop} \label{prop:hPlaneCrit}

For any $a \in \LL{n-3}$, $g$ is critical in the direction off of
the hyperplane.

\end{prop}

\begin{proof}

By symmetry, consider the \Lt{n-2} given by $\{u_1=0\}$.  For any $a
\in \{u_1=0\}$, the first row of $D\gT(a)$ vanishes. Thus, the local
behavior of $\widetilde{g}$ collapses points onto \Lt{n-2}. Explicit
calculation reveals that the collapse occurs in the direction of
$(2,1,\dots,1)$.

\end{proof}

Recall that the $p_m$ represent the point sets of \G{n} orbits
determined by intersecting  the \LL{n-3}.  Refer to these orbits
as ``$p_m$-points." First of all, each such point is
superattracting in all directions.

\begin{thm} \label{thm:pmSuperatt}

Under $g$, the fixed $p_m$-points are superattracting in every
direction.  Conversely, the only points that are superattracting
in every direction are the $p_m$-points.

\end{thm}

\begin{proof}

To establish that, at $p_m$, $g$ is critical in every direction in
\CP{n-1} show that the Jacobian $D\gT$ at $p_m$ has rank $1$.
Here, $p_m$ is lifted in the literal way. It then follows that,
since $\gT(p_m)\neq 0$, there are $n-2$ non-radial directions
through $p_m$ that have zero eigenvalue.

The Jacobian has the form
$$
D\gT =
\begin{pmatrix}
(a_{ij})&(b_{ij})\\
0&0
\end{pmatrix}
$$
where
\begin{align*}
a_{ij}=&\ \begin{cases}
3\,G_i(p_m)+\pdd{G_i}{u_i}(p_m)&i=j\\
\pdd{G_i}{u_j}(p_m)&i\neq j
\end{cases}
\quad i,j \leq m \\
b_{ij}=&\ \pdd{G_{i}}{u_{m+j}}(p_m) \quad i \leq m< j.
\end{align*}
With $S_k=S_{n,k}$ a straightforward calculation establishes that,
for $\ell\leq m$,
$$\pdd{S_k}{u_\ell}(p_m)=\begin{cases}
0&k>m\\
\binom{m-1}{k-1}&k\leq m
\end{cases}
$$
so that $\pdd{G_i}{u_j}(p_m)$ is the same value for $i,j\leq m$
with $i\neq j$.  Similarly, $\pdd{G_i}{u_\ell}(p_m)$ is the same
value for $\ell>m$.  It remains to show that
$$
3\,G_i(p_m)+\pdd{G_i}{u_i}(p_m)=\pdd{G_j}{u_k}(p_m) \quad
\text{for all}\ i,j,k\leq m.
$$

By manipulation of sums,
$$
\pdd{G_i}{u_i}(p_m)= \sum_{k=0}^{n-2} (-1)^{n-k}\,
\frac{n-k-1}{n-k+1}\,(n-2-k)S_k(p_m) + \sum_{k=0}^{n-2}
(-1)^{n-k}\,\frac{n-k-1}{n-k+1}\,\pdd{S_p}{u_i}(p_m).
$$
The second sum is $\pdd{G_j}{u_\ell}(p_m)$ for $j,\ell \leq m$ and
$j\neq \ell$. To show that the first sum amounts to
$-3\,G_i(p_m)$, notice that, from the proof of
Theorem~\ref{thm:gHlm},
\begin{align*}
\sum_{k=0}^{n-2} (-1)^{n-k}\, \frac{n-k-1}{n-k+1}\,(n-2-k)S_k(p_m)
&=\
(-1)^n \sum_{k=0}^{m}(-1)^k\,\frac{n-k-1}{n-k+2}\,(n-2-k)\binom{m}{k}\\
&=\ (-1)^n \sum_{k=0}^{m}
  (-1)^k\,\frac{n-k-1}{n-k+1}\,((n-k+1)-3)\binom{m}{k}\\
&=\ (-1)^n \sum_{k=0}^{m}
  (-1)^k\,(n-1-k)\,\binom{m}{k}-3\,G_i(p_m).
\end{align*}
Finally, the calculation reduces to showing that the first sum
vanishes.  This follows readily by splitting the sum into two
terms each of which is a binomial expansion of $1-1$.
Specifically,
\begin{align*}
\sum_{k=0}^{m} (-1)^k\,(n-1-k)\,\binom{m}{k} =&\
  (n-1)\sum_{k=0}^{m} (-1)^k\,\binom{m}{k}-
    \sum_{k=0}^{m} (-1)^k\,k\,\binom{m}{k}\\
    =&\ (n-1)(1-1)^m + m\sum_{k=1}^{m} (-1)^k\,\binom{m-1}{k-1}\\
    =&\ m\,(1-1)^{m-1}.
\end{align*}
Thus, the nonzero rows of $D\gT(p_m)$ are identical and the matrix
has rank $1$.

For the converse claim, consider a point $q$ that is critical in
every direction.  When $g$ is restricted to any intersection
\LL{k} of hyperplanes each of which is an \LL{n-2}, $q$ is again
critical for the restriction $g|_{\LL{k}}$.  Hence, $q$ lies on
some \LL{n-2} that does not contain \LL{k} and so, is determined
by the intersection of \LL{n-2} spaces.

\end{proof}

Now for the issue of the Fatou set $F_g$. Is there a Fatou
component of $g$ that is not in the basin of a $p_m$ point?

\begin{thm}

For $n=3,4$, $F_g$ consists of the basins of attraction of the
$p_m$-points.

\end{thm}

\begin{proof}

When $n=3$, the one-dimensional map $g$ has three fixed critical
$p_m$-points.  A basic result in one-dimensional dynamics states
that the Fatou set of a rational map with periodic critical points
consists only of superattracting basins; indeed, the basins have
full measure in \CP{1}.

In the two-dimensional case $n=4$, the claim follows from
Theorem~\ref{thm:pmSuperatt} and \cite{fs1}, Theorem~7.7. The latter
implies that if a holomorphic map $f$ on \CP{2} has a critical set
$C$ such that 1) $C$ is periodic and 2) $\CP{2}-C$ is Kobayashi
hyperbolic, then $f$ has only superattracting basins in its Fatou
set. See below for an explanation of the fact that condition 2)
applies to $g$.

\end{proof}

The general case remains open.

\begin{conj} \label{conj:gFatou}

For $n\geq 5$, $F_g$ consists of the basins of attraction of the
$p_m$-points.

\end{conj}

One approach to this claim adopts a technique from the proof of
Theorem~\ref{thm:gHlm}: reduction of dimension to the
one-dimensional case where some things are understood. The
argument for Theorem~\ref{thm:gFatouDense} employs the same idea.
Assume an arbitrary choice of $n \geq 5$.

The question of whether the basins of the $p_m$-points exhaust
$F_g$ calls for some preparation. Following \cite{ueda}, let $C_f$
be the critical set of a holomorphic map $f$ on \CP{m},
$$
D_f:=\bigcup_{k=1}^\infty f^k(C_f)\quad \text{and} \quad
E_f:=\bigcap_{k=1}^\infty f^k(\overline{D_f})
$$
be the \emph{postcritical set} and the $\omega$-limit set of $C_f$
respectively.  Also, the \emph{Fatou limit set} $\Lambda_f$ is
where the forward orbits of Fatou components accumulate. In the
case of $g$, $D_g=E_g=X$.

Let $p\in F_g$ and $U$ be the Fatou component to which $p$
belongs. For a critically-finite map $f$, $\Lambda_f \subset E_f$
(\cite{ueda}, Theorem 5.1). Accordingly, the forward orbit
$\{g^k(p)\}$ of $p$ accumulates on some \LL{n-3} and, by
Proposition~\ref{prop:hPlaneCrit}, is attracted to that
\LL{n-3}---call it \LL{n-3} as well.  Accordingly,
$$g^n(U) \longrightarrow \LL{n-3}.$$
The claim also follows from \cite{mcmull}, Theorem~2.36---a result
established by consideration of expansion in the Kobayashi metric
on the complement of the postcritical set.

The task now is to show that
$$g^r(U) \cap \LL{n-3} \neq \emptyset$$
for some $r$.  An argument might develop in two steps:  1) the
orbit of a point that is Fatou for $g$ accumulates at points that
are Fatou for $\gT:=g|_{\LL{n-3}}$; 2) a point that is Fatou for
\gT\ is also Fatou  for $g$ and, thereby, belongs to a Fatou
component in \CP{n-2}.

To treat the first claim, let $q\in \LL{n-3}$ be a limit point of
$\{g^k(p)\}$ with $g^{n_k}|_K \longrightarrow h$ where $h:K
\rightarrow \LL{n-3}$, $K\subset U$ is a neighborhood of $p$, and
$h(p)=q$.

Suppose that $q$ belongs to the Julia set $J_{\gT}$. By
Proposition~\ref{prop:hPlaneCrit} $g$ is superattracting at
$g^{k}(q)$ in some direction away from \LL{n-3} for all $k$. This
equips $q$ with a stable set
$$
S_q =
\{x\ |\ \mathrm{dist}(g^{n_k}(x),g^{n_k}(q)) \longrightarrow 0 \}
$$
transverse to \LL{n-3}. If \gT\ were hyperbolic---as in the case
$n=4$, one might expect that the Kobayashi expansion at $q$ would
produce saddle-like behavior and force $U$ to contain Julia points
for $g$.

To see claim 2) above, let $q \in F_{\gT}$ with a neighborhood
$\widetilde{N}$ on which $\{\gT^k\}$ is normal. Take $N$ to be the
connected neighborhood of $q$ that is absorbed by $\widetilde{N}$
and includes $\widetilde{N}$; that is, $N$ is the connected
component of the stable set of $\widetilde{N}$
$$S_{\widetilde{N}} = \bigcup_{x \in \widetilde{N}} S_x$$
where $\widetilde{N} \subset N$ and $S_x$ is the stable set of $x$.
Every point in $N$ belongs to some $S_x$. Thus, if $\gT^{n_k}$
converges to $\widetilde{h}$ on $\widetilde{N}$, then $g^{n_k}$
converges on $N$ to
$$h(y) = \widetilde{h}(x),\quad y \in S_x.$$

The claims 1) and 2) imply that some $g^r(U)$ intersects \LL{n-3};
indeed, $g^r(U) \cap \LL{n-3}$ is a Fatou component for \gT. By
the critical finiteness of \gT, the forward orbit of $g^r(U) \cap
\LL{n-3}$ meets some \LL{n-4} in Fatou points for $g|_{\LL{n-4}}$.

This cascade continues until some $g^s(U)$ makes contact with a
line \LL{1}, in particular, with the Fatou set of $g|_{\LL{1}}$.
Since $g|_{\LL{1}}$ has fixed critical points, it has only
superattracting basins. The only critical points on \LL{1} are
$p_m$-points. Hence, $g^s(U) \cap \LL{1}$ lies in the basin of
attraction of some such point.

How ``large" are the basins of the $p_m$-points?  First of all,
let
$$B_f := \bigcup_{k\geq 0} f^{-k}(C_f)$$
be the \emph{precritical} set of $f$. The following basic result
yields that the closure of $B_g$ contains the Julia set $J_g$.
(\cite{fs1}, Proposition 6.5.)

\begin{thm} \label{thm:Julia}

If $f:\CP{k}\rightarrow \CP{k}$ is holomorphic and
$\CP{k}-\overline{B_f}$ is hyperbolically embedded,

$$J_f \subset A_f:=\bigcap_{n>0}\ \overline{\bigcup_{m>n} f^{-m}(C_f)}.$$

\end{thm}
To apply this result to $g$, we must see that it satisfies the
hypotheses.  By Theorem~\ref{thm:gHlm}, $g$ is holomorphic on
\CP{n-2}. Two theorems of M. Green imply that
$\LL{n-1-m}-\overline{B_{g|_{\LL{n-1-m}}}}$ is hyperbolically
embedded in \LL{n-1-m} (taking $\LL{n-2}=\Hh$). (For details on
Green's results, consult \cite{fs1}, Section 5.) To see this,
suppose that, for $n \geq 4$ and $m\geq 2$,
$$\phi:\CC{} \longrightarrow \LL{n-1-m}-\overline{B_{g|_{\LL{n-1-m}}}}$$
is holomorphic.  Then $\phi(\CC{})$ omits at least $n-m+1$
hypersurfaces in \LL{n-1-m}, namely, some \LL{n-m-2} spaces and
their preimages. By one of Green's theorems (Theorem 5.6 in
\cite{fs1}), $\phi(\CC{})$ is contained in a compact complex
hypersurface. Since such a hypersurface intersects the omitted
hypersurfaces, $\phi(\CC{})$ omits at least three points and so,
is constant. The statement concerning hyperbolic embedding follows
from Green's other theorem (Theorem 5.5 in \cite{fs1}).

One other preliminary: since $C_g \subset g^{-1}(C_g),\ J_g
\subset A_g=\overline{B_g}$.  We can now establish a bit of
$F_g$'s global structure.

\begin{thm} \label{thm:gFatouDense}

Under the assumption that Conjecture~\ref{conj:gFatou} holds, the
Fatou set $F_g$ is dense in \Hh.

\end{thm}

\begin{proof}

Consider $j_0 \in J_g$ and let $U_0$ be a neighborhood of $j_0$.
By Theorem~\ref{thm:Julia}, some precritical points meet $U_0$ so
that, for some $m$,
$$g^m(U_0) \cap C_g \neq \emptyset.$$
If
$$U_1 := g^m(U_0) \cap \LL{n-3}$$
fails to contain Julia points, the case is made.  Otherwise, take
a Julia point $j_1 \in U_1$, a neighborhood of $j_1$.

The map $g|_{\LL{n-3}}$ is critically finite with critical set
$C_{n-3}$ in the intersection of \LL{n-3} and the hyperplanes in
$X$ different from \LL{n-3}. Hence, $C_{n-3}$ is a collection of
\LL{n-4} spaces. Implementing the argument given for $j_0$ and
$U_0$ under $g$ using $j_1$ and $U_1$ under $g|_{\LL{n-3}}$
produces a neighborhood of a Julia point $j_2$ on some \LL{n-4}.
The descent continues until it reaches a Julia point $j_{n-3}$ and
neighborhood $U_{n-3}$ on an \LL{1}.  Thus, $U_{n-3}$ meets the
Fatou set of $g|_{\LL{1}}$.  Since $g|_{\LL{1}}$ has fixed
critical points that are $p_m$-points, its Fatou set consists of
the superattracting basins of those $p_m$-points. Accordingly,
$U_{n-3}$---hence, $U_0$---contains points in $F_g$.

\end{proof}

\subsection{A query on the structure of g's Julia set}

For the restricted map $\widehat{g} = g|_{\LL{n-3}}$, the Julia
set is given by
$$J_{\widehat{g}} = J_g \cap \LL{n-3}.$$
The inclusion $J_{\widehat{g}} \subset J_g \cap \LL{n-3}$ is
clear. If $x \notin J_{\widehat{g}}$, then $x$ belongs to a basin
of a $p_m$-point so that $x \notin J_g$.  At each point $p\in
J_{\widehat{g}}$, the map is superattracting in the direction away
from \LL{n-3}. Thus, there is a ``stable set" $S_p$ of points in
$J_g$ whose orbits are attracted to the orbit of $p$. Accordingly,
there is a stable bundle over $J_{\widehat{g}}$
$$S_{J_{\widehat{g}}}:=\bigcup_{p \in J_{\widehat{g}}} S_p \subset J_g.$$
Are the $S_p$ one-dimensional manifolds? Are the preimages of the
$S_{J_{\widehat{g}}}$ dense in $J_g$?

In the case $n=4$, $g$ restricts to a critically-finite map
$\widehat{g}$ on an \LL{1} that is one of the six lines of
reflection for \G{4}. Figure~\ref{fig:g5CP1} displays the three
basins of attraction for $\widehat{g}$.  The Julia set
$J_{\widehat{g}}$ consists of the boundaries of these basins. For
each Julia point $p \in \LL{1}$, there is an $S_p$ away from the
line. What can be said about the structure of
$S_{J_{\widehat{g}}}$?

What about the points
$$K:=J_g-\bigcup_X S_{J_{\widehat{g}}}$$
that are not absorbed by $X$?

On an \LL{1}, each Julia point is non-wandering and has a
contracting direction onto \LL{1} and an expanding direction in
\LL{1}.  For a hyperbolic map on \CP{2}, the literature describes
a grading of the non-wandering set $\Omega$ by the expanding
dimension (\cite{fs2}):
$$\Omega=\Omega_0 \cup \Omega_1 \cup \Omega_2.$$
The $p_m$-points comprise $\Omega_0$ and $\cup_X J_{\widehat{g}}
\subset \Omega_1$. The non-wandering points not on $X$ belong to
$K$. Since any neighborhood of such a point $p$ contains an open
set that is attracted to $X$, there is expansion at $p$.  Does it
happen that
$$\Omega \cap K \subset \Omega_2$$
so that $g$ is hyperbolic?

\section{Geometry and dynamics in low-dimension}

To avoid confusion, let $g_{n+1}$ represent the particular map $g$
on the respective \G{n}-symmetric \Hh.

\subsection{The one-dimensional case: $g_4$ and Halley's method}

When $n=3$, the reflecting ``hyperplanes" consist of a three-point
orbit.  With these points located at
$$\{1,\rho,\rho^2 \ |\ \rho=e^{2\,\pi\,i/3}\},$$
the map's inhomogeneous expression on $\{u_2 \neq 0\}$ is
$$z \longrightarrow \frac{z(z^3- 2)}{2\,z^3 - 1}.$$

We can realize the \G{3} action on \CP{1} by the polyhedral
configuration of a double triangular pyramid---two regular
tetrahedra joined at a face.  The two-point orbit resides at $0$
and $\infty$ and defines two hemispheres in the usual way.
Accordingly, the unit circle corresponds to the equatorial
boundary between hemispheres and the $3$-points
$\{1,\rho,\rho^2\}$ are vertices where four faces congregate.

\begin{figure}[h]

\resizebox{\basinWidth}{!}{\includegraphics{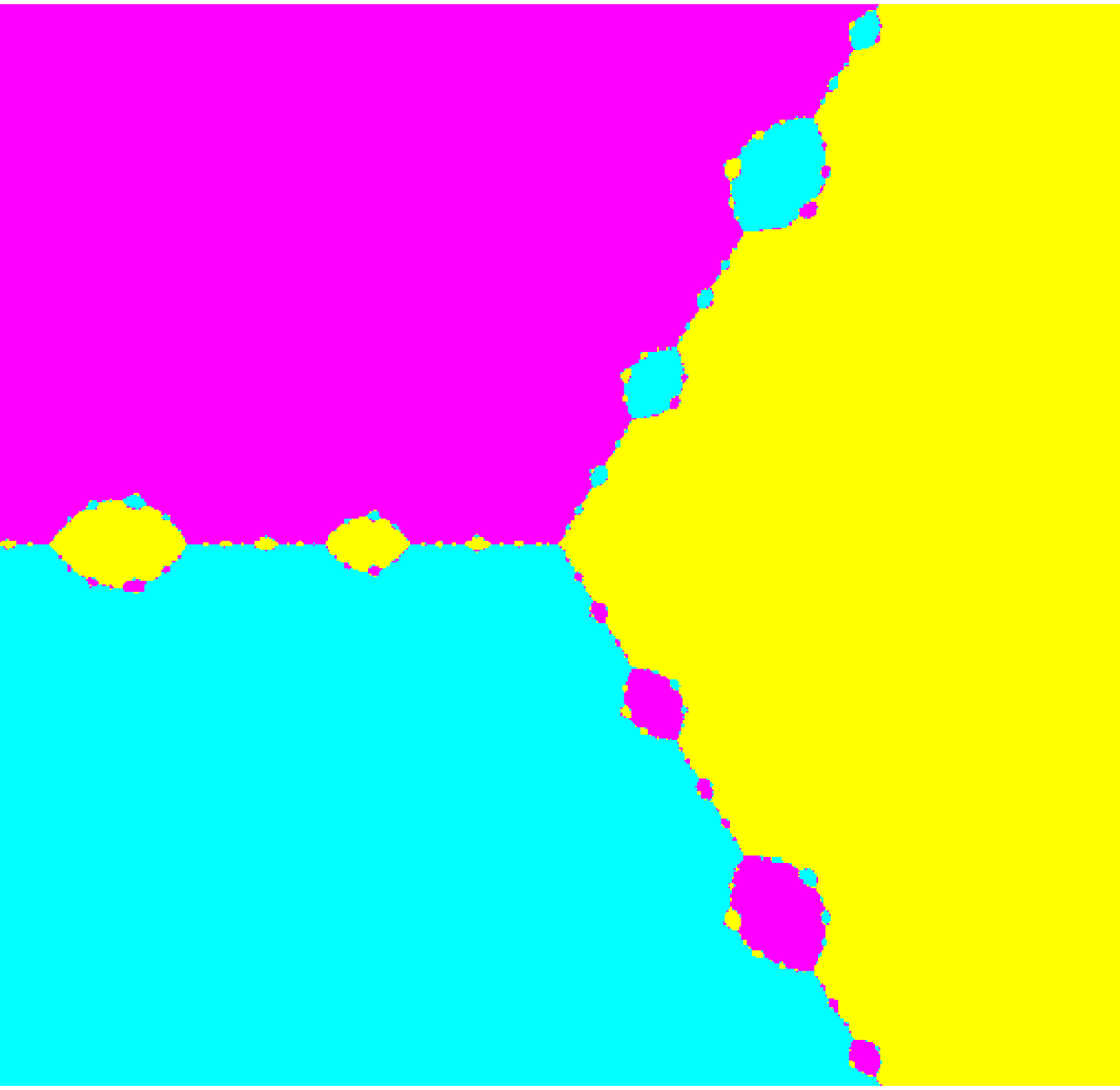}}

\caption{Dynamics of $g_4$ on the \Sym{3}-symmetric \CP{1}}

\label{fig:g4}

\end{figure}

Consider the degree-$4$ map that fixes the vertices of each face
and sends one face \textsf{F} to four others: \textsf{F} itself
and the three faces in the hemisphere not containing \textsf{F}.
This symmetrical construction results in \G{3}-equivariant
behavior. At the three equatorial vertices, the map opens up a
face's internal angle of $\pi/2$ to an angle of $3\,\pi/2$ so that
the local behavior is cubing.  This makes the $3$-point orbit
doubly-critical and, by degree counting, the entire critical set.
Accordingly, this map must be $g_4$.  Since $g_4$ has periodic
critical points, the superattracting basins constitute its Fatou
set and, moreover, have full measure in \CP{1}.  A portrait of the
basins appears in Figure~\ref{fig:g4}.

It turns out that $g_4$ is Halley's Method---a variation on
Newton's Method---for a cubic polynomial.   (See \cite{halley} for
a description of Halley's Method in real variables.)  In the
coordinates selected above, the polynomial to which we apply
Halley's method is
$$z^3-1.$$

%%%%%%%%%%%%%%
%%%%%%%%%%%%%%

\subsection{The map in two dimensions}

Since $g_{n+1}$ has real coefficients, it preserves the \RP{n-2}
of points whose coordinates can be expressed by real numbers. Call
this space $R$.  Under \G{4}, $R$ has the structure of a
projective cube.  We can view this as a hemisphere where one
vertex is at the pole and the other three vertices lie along a
circle whose center is the distinguished vertex.  The $3$-point
orbit (i.e., the face-centers) lies on another circle centered at
the north pole.

\begin{figure}[h]

\resizebox{\basinWidth}{!}{\includegraphics{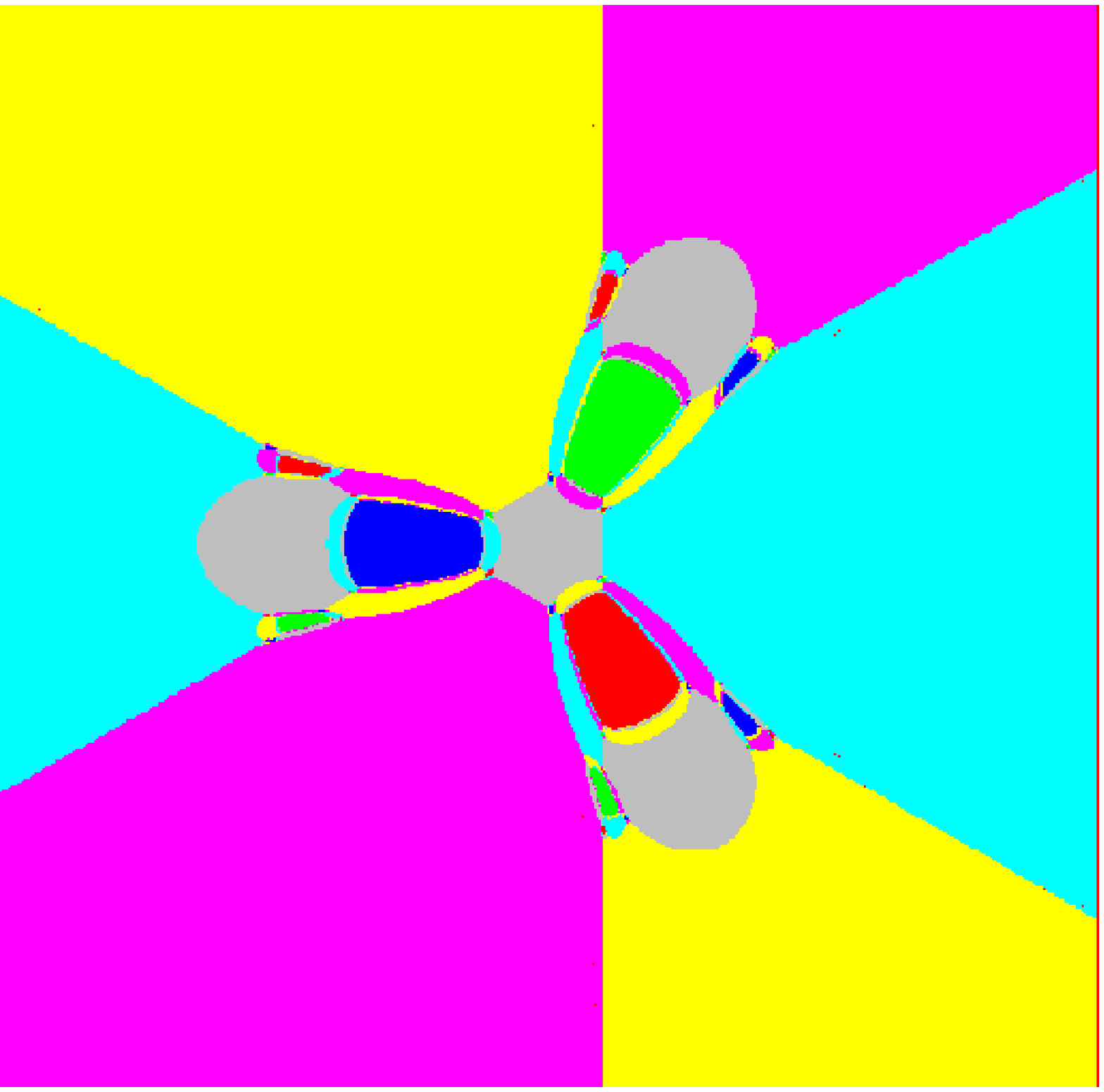}}

\caption{Dynamics of $g_5$ on the \Sym{4}-symmetric \RP{2}}

\label{fig:g5RP2}

\end{figure}

Figure~\ref{fig:g5RP2} displays the basins of attraction of $g_5$
on $R$.  In the affine plane of the picture,  the vertices of the
cube are
$$(0,0), (1,0), \Biggl(-\frac{1}{2},\pm \frac{\sqrt{3}}{2}\Biggr)$$
while the three face-centers are the edge-midpoints $$
\Bigl(-\frac{1}{2},0\Bigr), \Biggl(\frac{1}{4}, \pm
\frac{\sqrt{3}}{4}\Biggr)
$$
of the equilateral triangle formed by the three vertices that are
not $(0,0)$.  The map is given by
\tiny
\begin{align*} (x,y)
\longrightarrow & \Biggl( \frac{3\,\bigl( 15\,x^4 + 12\,x^5 -
       30\,x^2\,y^2 - 5\,y^4 +
       20\,{x}\,y^4 \bigr) }{1 - 10\,x^2 +
     20\,x^3 + 30\,x^4 + 40\,x^5 -
     10\,y^2 - 60\,{x}\,y^2 +
     60\,x^2\,y^2 - 80\,x^3\,y^2 +
     30\,y^4 - 120\,{x}\,y^4},\\
&  \frac{24\,\bigl( -5\,{x}\,y^3 +
       5\,x^2\,y^3 + y^5 \bigr) }{1 -
     10\,x^2 + 20\,x^3 + 30\,x^4 +
     40\,x^5 - 10\,y^2 - 60\,{x}\,y^2 +
     60\,x^2\,y^2 - 80\,x^3\,y^2 +
     30\,y^4 - 120\,{x}\,y^4} \Biggr).
\end{align*}
\normalsize

\begin{figure}[h]

\resizebox{2in}{!}{\includegraphics{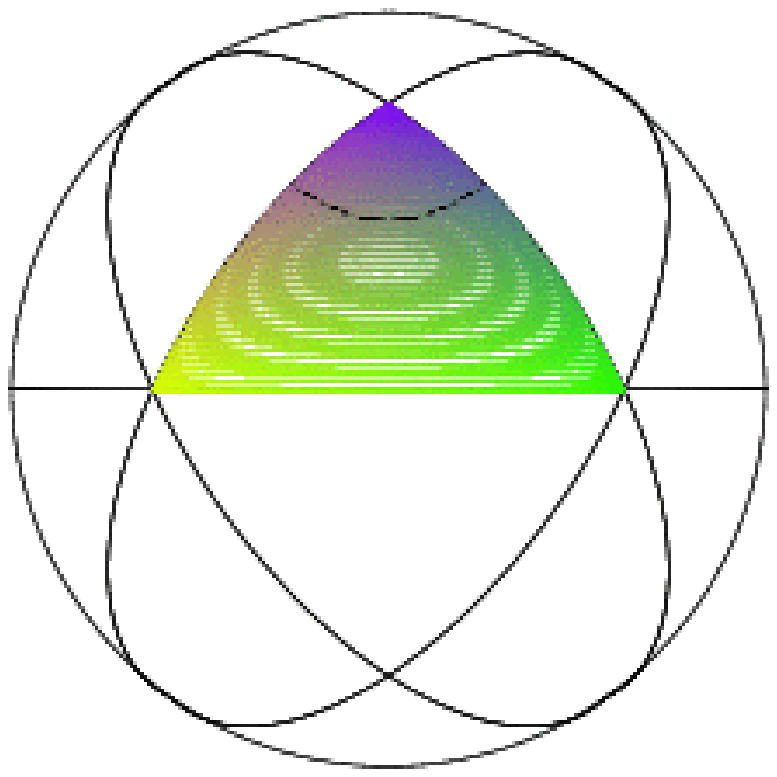}}
$\longrightarrow$
\resizebox{2in}{!}{\includegraphics{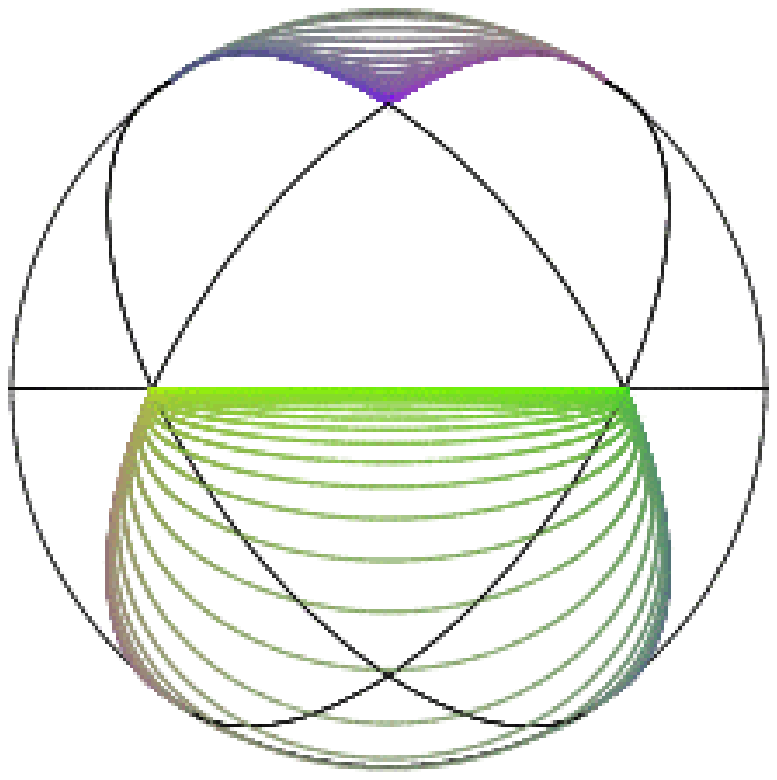}}

\caption{Geometry of $g_5$ on the \Sym{4}-symmetric \RP{2}}

\label{fig:g5Geom}

\end{figure}

The six lines of reflection run along the edges and a diagonal of
a face.  These lines carve the hemisphere into twelve triangles
each of which is a fundamental domain for the reflection group
action \G{4}.  Viewing the ``hemi-cube" from above an edge,
Figure~\ref{fig:g5Geom} reveals the map's action on a fundamental
triangle: one triangle stretches and twists onto five other
associated triangles.

%%%%%%%%%%%%%%
%%%%%%%%%%%%%%
%%%%%%%%%%%%%%

\begin{figure}[h]

\resizebox{\basinWidth}{!}{\includegraphics{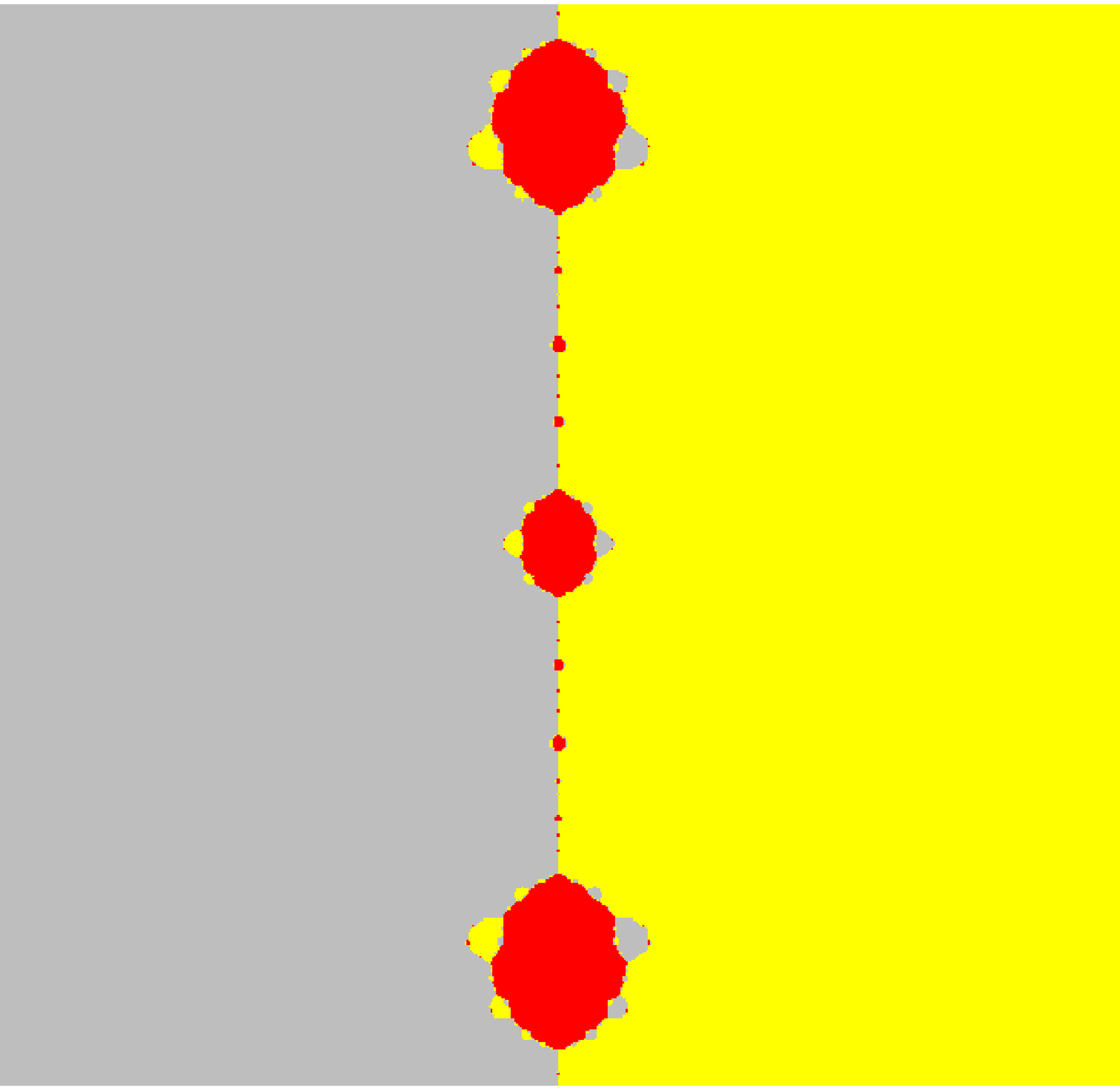}}

\caption{Dynamics of $g_5$ on the \Z{2}-symmetric \CP{1}}

\label{fig:g5CP1}

\end{figure}

Returning to $u$ coordinates, one of the six mirrors---say,
$\{u_3=0\}$---is \Z{2}-stable. Restricted to this line, $g_5$ has
three superattracting points:

\begin{itemize}

\item a two-point \Z{2} orbit of type $p_1$ points $[1,0,0]$ and
$[0,1,0]$\\
(where $\{u_2=0\}$ and $\{u_1=0\}$ intersect $\{u_3=0\}$)

\item a one-point \Z{2} orbit of  the point $p_2=[1,1,0]$\\
(where $\{u_1=u_2\}$ intersects $\{u_3=0\}$).

\end{itemize}
In coordinates where the two-point orbit is $\pm 1$ and the
one-point orbit is $0$, the map takes the form
$$z \longrightarrow \frac{4\,z^3(z^2+5)}{15\,z^4+10\,z^2-1}.$$
Figure~\ref{fig:g5CP1} shows their basins of attraction on the
line.  Notice that this \CP{1} intersects $R$ in an \RP{1} that
corresponds to a line of reflective symmetry in
Figure~\ref{fig:g5RP2} and the horizontal mirror in
Figure~\ref{fig:g5CP1}---for instance, the line that passes
through the red, gray, and yellow basins.
%%%%%%%%%%%%%%
%%%%%%%%%%%%%%

\subsection{The three-dimensional map: A cascade of critical
finiteness}

\begin{figure}[h]

\resizebox{\basinWidth}{!}{\includegraphics{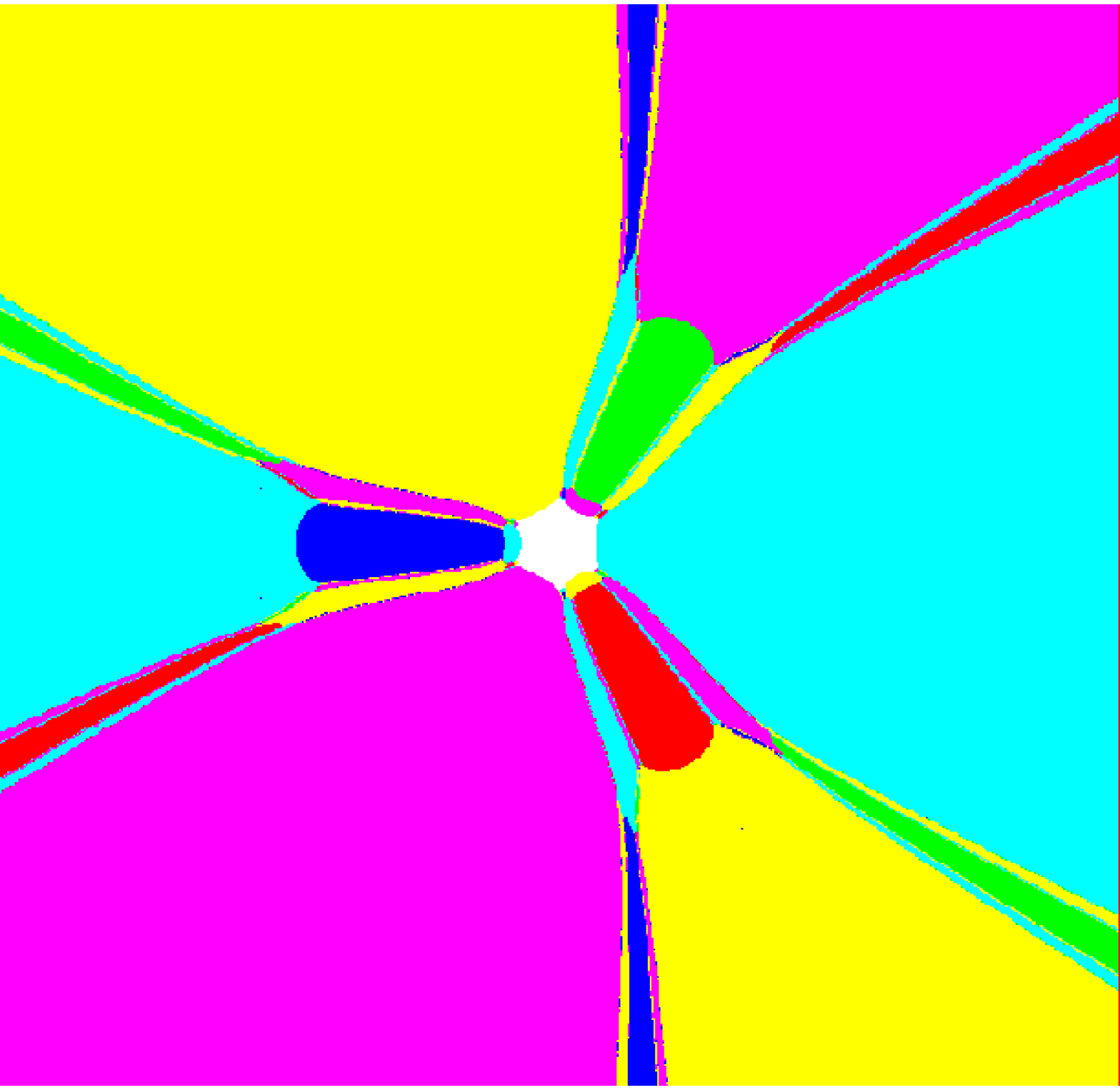}}

\caption{Dynamics of $g_6$ on the \Sym{3}-symmetric \RP{2}}

\label{fig:g6RP2}

\end{figure}

A component of $g_6$'s critical set is a \CP{2}.  On the
\Sym{3}-symmetric $\{u_4=0\}$ the map has three \Sym{3} orbits of
superattracting points:

\begin{itemize}

\item type $p_1$ points $[1,0,0,0],[0,1,0,0],[0,0,1,0]$

\item type $p_2$ points $[1,1,0,0],[1,0,1,0],[0,1,1,0]$

\item $p_3=[1,1,1,0]$.

\end{itemize}
In the basin plot on the corresponding \RP{2}
(Figure~\ref{fig:g6RP2}), the geometry is that of a projective
double triangular pyramid and these points respectively occupy
$$
(1,0), \Biggl(-\frac{1}{2},\pm \frac{\sqrt{3}}{2}\Biggr) \qquad
\Bigl(-\frac{1}{2},0\Bigr), \Biggl(\frac{1}{4}, \pm
\frac{\sqrt{3}}{4}\Biggr) \qquad (0,0).
$$  The map is given by
\tiny
\begin{align*}
(x,y) &\longrightarrow \\
& \biggl( 9\,\Bigl( 15\,x^4 + 24\,x^5 +
       15\,x^6 - 30\,x^2\,y^2 -
       15\,x^4\,y^2 - 5\,y^4 +
       40\,x\,y^4 - 35\,x^2\,y^4 -
       5\,y^6 \Bigr),\\
& -72\,y^3\,\Bigl( 5\,x - 10\,x^2 +
       5\,x^3 - 2\,y^2 + 5\,x\,y^2
       \Bigr)\biggr)/ \\
&\Bigl(1 - 15\,x^2 + 40\,x^3 +
     90\,x^4 + 240\,x^5 + 130\,x^6 -
     15\,y^2 - 120\,x\,y^2 +
     180\,x^2\,y^2\\
&- 480\,x^3\,y^2 + 30\,x^4\,y^2  + 90\,y^4 - 720\,x\,y^4 +
     630\,x^2\,y^4 + 90\,y^6 \Bigr) .
\end{align*}
\normalsize This image makes for interesting comparison to the
\Sym{4}-symmetric Figure~\ref{fig:g5RP2}.
%%%%%%%%%%%%%%
%%%%%%%%%%%%%%

\begin{figure}[h]

\resizebox{\basinWidth}{!}{\includegraphics{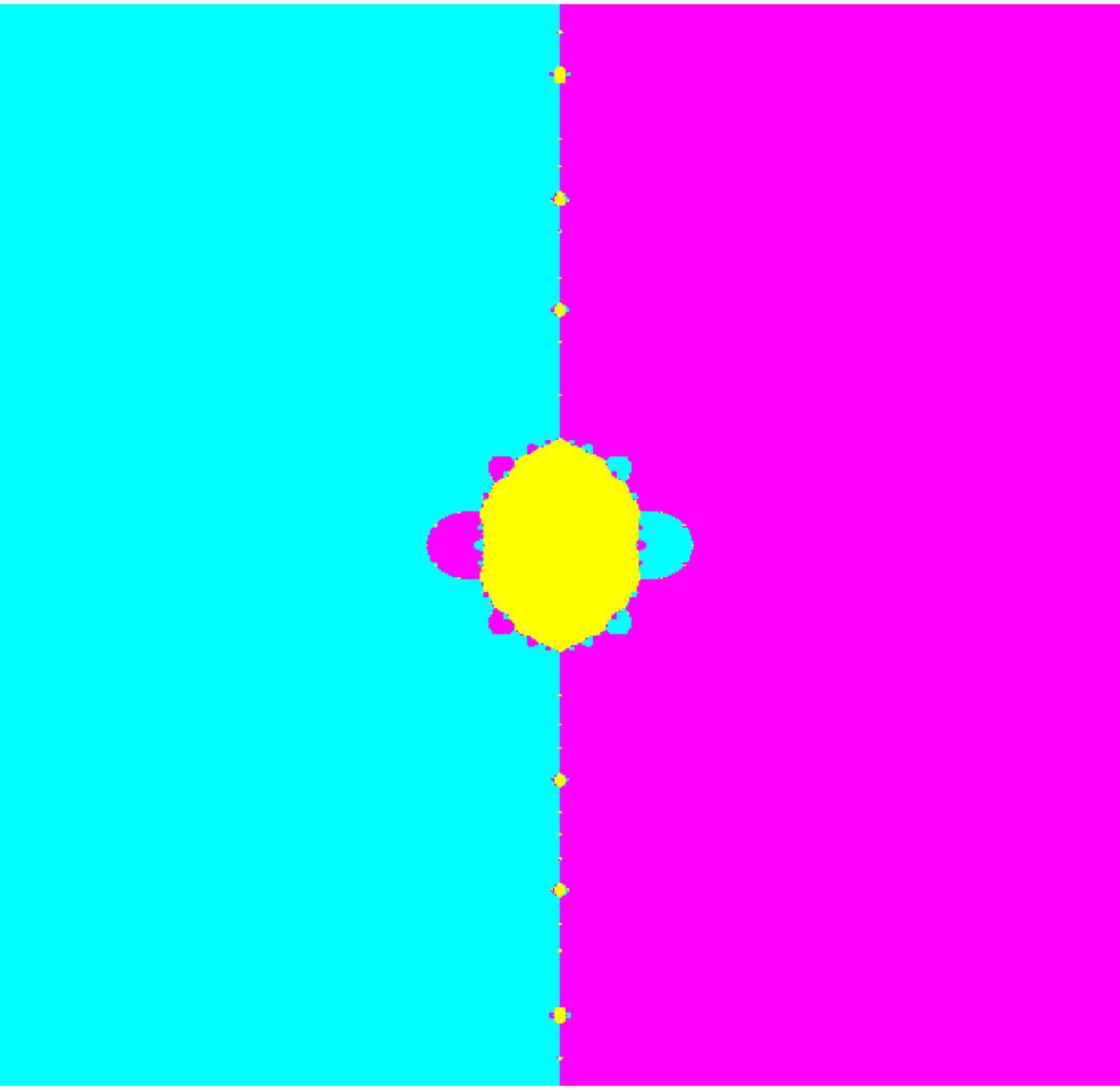}}

\caption{Dynamics of $g_6$ on the \Z{2}-symmetric \CP{1}}

\label{fig:g6CP1Z2}

\end{figure}

\begin{figure}[h]

\resizebox{\basinWidth}{!}{\includegraphics{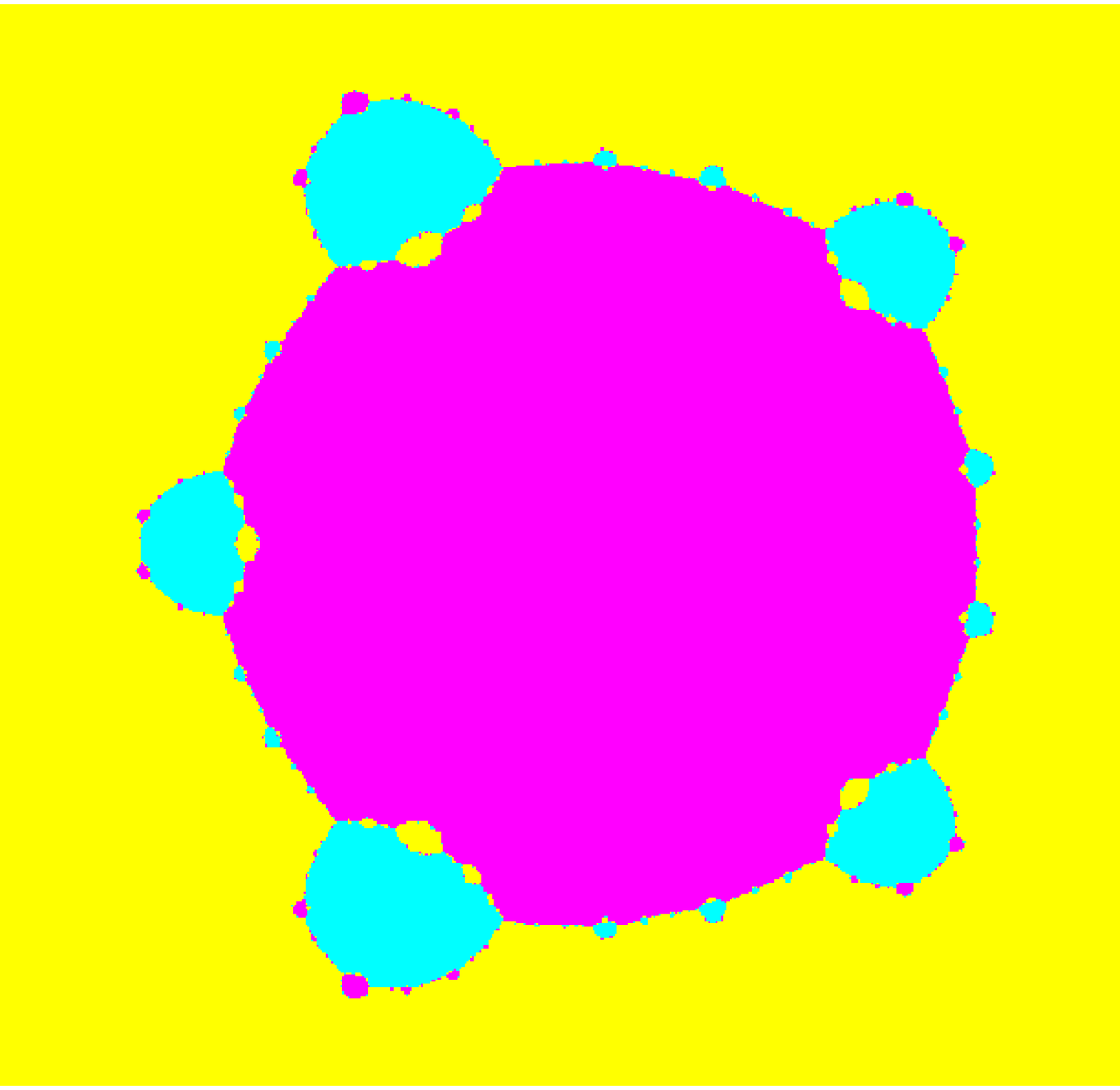}}

\caption{Dynamics of $g_6$ on the \Z{1}-symmetric \CP{1}}

\label{fig:g6CP1Z1}

\end{figure}

On the critical component $\{u_4=0\}$, $g|_{\{u_4=0\}}$ has two
types of critical line: $\{u_3=0\}$ and $\{u_2=u_3\}.$ The
respective lines have \Z{2} and trivial symmetry.  As for
superattracting points, the former line contains
$[1,0,0,0],[0,1,0,0]$ (a two-point \Z{2} orbit) and $[1,1,0,0]$
while on the latter line we find $[1,0,0,0],[0,1,1,0],[1,1,1,0]$.
In the respective basin plots for $g_6$ restricted to the lines
(Figure~\ref{fig:g6CP1Z2} and Figure~\ref{fig:g6CP1Z1}), these
points are $\pm 1,0,$ and $1,0,-1$ while the maps are
$$
z \longrightarrow
\frac{8\,z^3(3\,z^2+5)}{5\,z^6+45\,z^4+15\,z^2-1}
$$
and
$$
z \longrightarrow
\frac{8\,z^4(z^2-2\,z+5)}{5\,z^6+30\,z^5+15\,z^4+20\,z^3-5\,z^2-2\,z+1}.
$$
As before, each \CP{1} intersects the \RP{2} of
Figure~\ref{fig:g6RP2} in an \RP{1}: the three lines
$$\{u_k=0\ |\ k=1,2,3\}$$
give the edges of the ``triangle" whose vertices are
$$(1,0),\Biggl(-\frac{1}{2},\pm \frac{\sqrt{3}}{2}\Biggr)$$
and the three lines
$$\{u_k=u_\ell\ |\ k,\ell=1,2,3\}$$
correspond to the lines of reflective symmetry through $(0,0)$.

%%%%%%%%%%%%%%
%%%%%%%%%%%%%%

%%%%%%%%%%%%%%
%%%%%%%%%%%%%%

\clearpage
\appendix
\section*{Appendix: Proofs of computational statements}

\begin{appLemma}[\textbf{\ref{lm:Snk}}]

For $n\geq 3$ and $k\leq n$, the \G{n-1} invariants $S_{n,k}$
transform under $T$ according to
$$
S_{n,k}(T u) = \sum_{\ell=0}^k (-1)^\ell \binom{n-k+\ell}{n-k}
u_1^\ell\,S_{n,k-\ell}.
$$

\end{appLemma}

\begin{proof}

The argument is induction on $n$.  Note first that
$$S_{n,0}=1$$
satisfies the identity trivially for all $n$.  At the other
extreme,
$$ S_{n,n}=0$$
also satisfies the statement.  To see this, examine
\begin{equation*}
 \sum_{\ell=0}^{n}
 (-1)^\ell\,\binom{\ell}{0}\,u_1^\ell\,S_{n,n-\ell}
 = \sum_{\ell=0}^{n}
 (-1)^\ell\,u_1^\ell\,S_{n,n-\ell}
 = -u_1 \sum_{m=0}^{n-1}
 (-1)^m\,u_1^m\,S_{n,n-1-m}.
\end{equation*}
For the final equality, use $S_{n,n}=0$ and set $\ell=m+1$.  By
substituting $x$ for the variable $u_1$ that appears explicitly,
the sum factors:
$$
\sum_{m=0}^{n-1} (-1)^m\,x^m\,S_{n,n-1-m} = \prod_{k=1}^{n-1}
(u_k-x).
$$
Consequently, it vanishes when $x=u_1$.

For the base $n=3$,
\begin{align*}
S_{3,1}(T u) &=\ -u_1 + (u_2-u_1)\\
&=\ (u_1 + u_2) - 3\,u_1\\
&=\ S_{3,1} - 3\,u_1
\end{align*}
and
\begin{align*}
S_{3,2}(T u) &=\ (-u_1)(u_2-u_1)\\
&=\ u_1 u_2 - 2\,(u_1^2 + u_1 u_2)  + 3\,u_1^2\\
&=\ S_{3,2} - 2\,u_1 S_{3,1} + 3\,u_1^2.
\end{align*}

To make the inductive step, use the reduction
$$S_{n+1,k} = S_{n,k} + u_n\,S_{n,k-1}$$
and assume the claim holds for $S_{n,k}$ and $S_{n,k-1}$. (Note
that the cases $k=n$ and $k=1$ fall under the scope of the remarks
above.) Thus,
\begin{align*}
S_{n+1,k}(T u) =&\ S_{n,k}(T u) + (u_n-u_1)\,S_{n,k-1}(T u) \\
=&\ \sum_{\ell=0}^k (-1)^\ell \binom{n-k+\ell}{n-k}
  u_1^\ell\,S_{n,k-\ell} \\
&\ - \sum_{m=0}^{k-1} (-1)^m \binom{n-(k-1)+m}{n-(k-1)}
  u_1^{m+1}\,S_{n,k-1-m} \\
&\ + \sum_{p=0}^{k-1} (-1)^p \binom{n-(k-1)+p}{n-(k-1)}
  u_1^p\,u_n\,S_{n,k-1-p} \\
=&\ S_{n,k} + \sum_{\ell=1}^k (-1)^\ell
\binom{n-k+\ell}{n-k}u_1^\ell\,S_{n,k-\ell} \\
&\ + \sum_{m=0}^{k-1} (-1)^{m+1} \binom{n-k+(m+1)}{n-k+1}
  u_1^{m+1}\,S_{n,k-(m+1)} \\
&\ + u_n S_{n,k-1} + \sum_{p=1}^{k-1} (-1)^p
\binom{n+1-k+p}{n+1-k} u_1^p\,u_n\,S_{n,k-p-1}.
\end{align*}
Setting $m=\ell-1$ and $p=\ell$ gives
\begin{align*}
S_{n+1,k}(T u) =&\ S_{n,k}+ u_n\,S_{n,k-1}\\
&\ + \sum_{\ell=1}^k (-1)^\ell
\binom{n-k+\ell}{n-k}u_1^\ell\,S_{n,k-\ell} \\
&\ + \sum_{\ell=1}^k (-1)^\ell \binom{n-k+\ell}{n-k+1}
  u_1^{\ell}\,S_{n,k-\ell} \\
&\  + \sum_{\ell=1}^{k-1} (-1)^\ell
\binom{n+1-k+\ell}{n+1-k}u_1^\ell\,u_n\,S_{n,k-\ell-1}\\
=&\ S_{n+1,k} + \sum_{\ell=1}^k (-1)^\ell
\Biggl(\binom{n-k+\ell}{n-k}+\binom{n-k+\ell}{n-k+1}\Biggr)
  u_1^\ell\,S_{n,k-\ell} \\
&\ + \sum_{\ell=1}^{k-1} (-1)^\ell
\binom{n+1-k+\ell}{n+1-k}u_1^\ell\,u_n\,S_{n,k-\ell-1}\\
=&\ S_{n+1,k} + \sum_{\ell=1}^k (-1)^\ell
\binom{n+1-k+\ell}{n+1-k}u_1^\ell\,S_{n,k-\ell} \\
&\ + \sum_{\ell=1}^{k-1} (-1)^\ell
\binom{n+1-k+\ell}{n+1-k}u_1^\ell\,u_n\,S_{n,k-\ell-1}\\
=&\ S_{n+1,k} + \sum_{\ell=1}^{k-1} (-1)^\ell
\binom{n+1-k+\ell}{n+1-k}u_1^\ell\,
  \Bigl(S_{n,k-\ell}+u_n\,S_{n,k-\ell-1}\Bigr) \\
& + (-1)^k\,\binom{n+1}{n+1-k}u_1^k  \\
=&\ S_{n+1,k} + \sum_{\ell=1}^{k-1} (-1)^\ell
\binom{n+1-k+\ell}{n+1-k}u_1^\ell\,S_{n+1,k-\ell}  +
(-1)^k\,\binom{n+1}{n+1-k}\,u_1^k  \\
=&\ \sum_{\ell=0}^{k} (-1)^\ell
\binom{(n+1)-k+\ell}{(n+1)-k}u_1^\ell\,S_{n+1,k-\ell}.
\end{align*}

\end{proof}

\begin{appLemma}[\textbf{\ref{lm:specSum}}]

$$
\sum_{k=0}^m (-1)^k\,\frac{k+1}{(k+3)!\,(m-k)!} =
\frac{m+1}{(m+3)!}.
$$

\end{appLemma}

\begin{proof}

Consider the expansion of the generating function
\begin{align*}
\frac{(1-x)^{m+3}}{x^2}
=&\ \sum_{\ell=0}^{m+3} (-1)^\ell\,\binom{m+3}{\ell}\,x^{\ell-2}\\
=&\ x^{-2} - (m+3)\,x^{-1} + \binom{m+3}{2} +
  \sum_{\ell=3}^{m+3} (-1)^{\ell-2}\,\binom{m+3}{\ell}\,x^{\ell-2}\\
=&\ x^{-2} - (m+3)\,x^{-1} + \binom{m+3}{2} +
  \sum_{k=0}^m
  (-1)^{k+1}\,\binom{m+3}{k+3}\,x^{k+1}.
\end{align*}
Now, differentiate and evaluate at $x=1$:
\begin{align*}
\dd{x}\biggl(\frac{(1-x)^{m+3}}{x^2}\biggr)\biggr|_{x=1}
=&\ -2 + m+3 + \sum_{k=0}^m
  (-1)^{k+1}\,(k+1)\binom{m+3}{k+3}\\
0=&\ m + 1 - (m+3)!\,\sum_{k=0}^m
  (-1)^k\,\frac{k+1}{(k+3)!\,(m-k)!}.
\end{align*}
Rearranging this equation yields the desired statement.

\end{proof}

\begin{appLemma}[\textbf{\ref{lm:g2specSum}}]

\begin{equation*}
\sum_{k=0}^m \frac{k+1}{k+3}\,\binom{m+2}{k+2}\,
  u_1^{m-k}\,(u_2-u_1)^{k+3}
= \frac{m+1}{m+3}\,\bigl(u_2^{m+3}-u_1^{m+3}\bigr)-
  u_1\,u_2\,\bigl(u_2^{m+1}-u_1^{m+1}\bigr).
\end{equation*}

\end{appLemma}

\begin{proof}

Letting $u=u_1$ and $v=u_2-u_1$,

\begin{align*}
\sum_{k=0}^m& \frac{k+1}{k+3}\,\binom{m+2}{k+2}\,
  u_1^{m-k}\,(u_2-u_1)^{k+3} = \frac{1}{m+3} \sum_{k=0}^m
  (k+1)\,\binom{m+3}{k+3}\,u^{m-k}\,v^{k+3} \\
=&\ \frac{1}{m+3}\,\Biggl(
    \sum_{k=0}^m (k+4)\,\binom{m+3}{k+3}\,
  u^{(m+3)-(k+3)}\,v^{k+3}
- 3\,\sum_{k=0}^m \binom{m+3}{k+3}\,
  u^{(m+3)-(k+3)}\,v^{k+3} \Biggr) \\
=&\ \frac{1}{m+3}\,\Biggl(
    \sum_{p=3}^{m+3} (p+1)\,\binom{m+3}{p}\,
  u^{(m+3)-p}\,v^p
- 3\,\sum_{p=3}^{m+3} \binom{m+3}{p}\,
  u^{(m+3)-p}\,v^p \Biggr) \\
=&\ \frac{1}{m+3}\,\Biggl(
    \sum_{p=0}^{m+3} (p+1)\,\binom{m+3}{p}\,
  u^{(m+3)-p}\,v^p
- 3\,\sum_{p=0}^{m+3} \binom{m+3}{p}\,
  u^{(m+3)-p}\,v^p \\
&-u^{m+3} - 2\,(m+3)\,u^{m+2}\,v -
3\,\binom{m+3}{2}\,u^{m+1}\,v^2\\ &+ 3\,u^{m+3} +
3\,(m+3)\,u^{m+2}\,v + 3\,\binom{m+3}{2}\,u^{m+1}\,v^2 \Biggr)
\\
=&\ \frac{1}{m+3}\,\Biggl(
    \sum_{p=0}^{m+3} (p+1)\,\binom{m+3}{p}\,
  u^{(m+3)-p}\,v^p
- 3\,\sum_{p=0}^{m+3} \binom{m+3}{p}\,
  u^{(m+3)-p}\,v^p \\
&+ 2\,u^{m+3} + (m+3)\,u^{m+2}\,v \Biggr).
\end{align*}
The second sum amounts to the binomial expansion of
$(u+v)^{m+3}=u_2^{m+3}$ while the first sum is the $v$-derivative
of $(u+v)^{m+3}\,v$.  In explicit terms, note that
$$
\sum_{p=0}^{m+3} \binom{m+3}{p}\,
  u^{(m+3)-p}\,v^{p+1} = (u+v)^{m+3}\,v.
$$
Hence,
\begin{align*}
\frac{\partial}{\partial v} \Biggl(
  \sum_{p=0}^{m+3} \binom{m+3}{p}\,
  u^{(m+3)-p}\,v^{p+1} \Biggr)
=&\ \frac{\partial}{\partial v} ((u+v)^{m+3}\,v)\\
\sum_{p=0}^{m+3} (p+1)\,\binom{m+3}{p}\,
  u^{(m+3)-p}\,v^p
=&\ (u+v)^{m+3} + (m+3)\,(u+v)^{m+2}\,v.
\end{align*}
Substituting into the expression above and reverting to $u_1$ and
$u_2$,
\begin{align*}
\sum_{k=0}^m& \frac{k+1}{k+3}\,\binom{m+2}{k+2}\,
  u_1^{m-k}\,(u_2-u_1)^{k+3} \\
=&\ \frac{1}{m+3} \Bigl( u_2^{m+3} + (m+3) \,u_2^{m+2}\,(u_2-u_1)
- 3\,u_2^{m+3} + 2\,u_1^{m+3}
+(m+3)\,u_1^{m+2}\,(u_2-u_1) \Bigr) \\
=&\ \frac{m+1}{m+3}\,\bigl(u_2^{m+3} - u_1^{m+3}\bigr)
  -u_1\,u_2\,\bigl(u_2^{m+1} - u_1^{m+1}\bigr).
\end{align*}

\end{proof}

\begin{appLemma}[\textbf{\ref{lm:G1pm}}]

$$
\sum_{p=0}^{m} (-1)^p\,\frac{n-p-1}{n-p+1}\, \binom{m}{p}
=\frac{2\,(-1)^{m-1}}{(n+1)\binom{n}{m}}.
$$

\end{appLemma}

\begin{proof}

Let
$$
\Lambda_{n,m}=\sum_{p=0}^{m} (-1)^p\,\frac{n-p-1}{n-p+1}\,
\binom{m}{p} \quad \text{and} \quad
L_{n,m}=\frac{2\,(-1)^{m-1}}{(n+1)\binom{n}{m}}.
$$

From the reduction
$$L_{n,m}=L_{n,m-1}-L_{n-1,m-1},$$
proceed by induction on $n$ and $m$.  For the base relative to
$m$:
$$
\Lambda_{n,1} = \frac{n-1}{n+1}-\frac{n-2}{n} = \frac{2}{n\,(n+1)}
= L_{n,1}.
$$
Make the inductive step by verifying that $\Lambda_{n,m}$ admits
the same reduction as $L_{n,m}$.  Consider
\begin{align*}
\Lambda_{n,m-1}-\Lambda_{n-1,m-1}=&\ \sum_{p=0}^{m-1}
(-1)^p\,\frac{n-p-1}{n-p+1}\,
\binom{m-1}{p} \\
&\ - \sum_{p=0}^{m-1} (-1)^p\,\frac{n-p-2}{n-p}\, \binom{m-1}{p}.
\end{align*}
Shearing the second sum by one term,
\begin{align*}
\Lambda_{n,m-1}-\Lambda_{n-1,m-1}=&\ \sum_{p=0}^{m-1}
(-1)^p\,\frac{n-p-1}{n-p+1}\,
\binom{m-1}{p} \\
&\ - \sum_{p=1}^m (-1)^{p-1}\,\frac{n-(p-1)-2}{n-(p-1)}\,
\binom{m-1}{p-1} \\
=&\ \frac{n-1}{n+1} + \sum_{p=1}^{m-1}
(-1)^p\,\frac{n-p-1}{n-p+1}\,
\binom{m}{p} \frac{m-p}{m} \\
&\ - \sum_{p=1}^{m-1} (-1)^{p-1}\,\frac{n-p-1}{n-p+1}\,
\binom{m}{p}\frac{p}{m} + (-1)^m\,\frac{n-m-1}{n-m+1} \\
=&\ \frac{n-1}{n+1} + \sum_{p=1}^{m-1} (-1)^p\,\frac{n-p-1}{n-p+1}
\binom{m}{p}
+ (-1)^m\,\frac{n-m-1}{n-m+1} \\
=&\ \sum_{p=1}^m (-1)^p\,\frac{n-p-1}{n-p+1} \binom{m}{p}\\
=&\ \Lambda_{n,m}.
\end{align*}

\end{proof}

\end{document}